\documentclass[a4paper,fleqn]{cas-sc}

\usepackage[authoryear,longnamesfirst]{natbib}
\usepackage{float}
\usepackage{placeins}
\usepackage{caption}
\usepackage{subcaption}
\RequirePackage[mathscr]{eucal}
\usepackage{longtable,tabularx}
\usepackage{nomencl}
\usepackage{colortbl}
\usepackage{amsmath}

\usepackage{atbegshi}
\AtBeginDocument{\AtBeginShipoutNext{\AtBeginShipoutDiscard}}

\def\tsc#1{\csdef{#1}{\textsc{\lowercase{#1}}\xspace}}
\tsc{WGM}
\tsc{QE}
\tsc{EP}
\tsc{PMS}
\tsc{BEC}
\tsc{DE}

\begin{document}
\let\WriteBookmarks\relax
\def\floatpagepagefraction{1}
\def\textpagefraction{.001}
\shorttitle{ }
\shortauthors{Ahmed et~al.}

\title [mode = title]{Impact of Separation Distance on the Performance and Annual Energy Production of a Dual-Flap Oscillating Surge Wave Energy Converter}

\author[1]{ Alaa Ahmed}

\cormark[1]

\affiliation[1]{Department of Civil, Environmental and Ocean Engineering, Stevens Institute of Technology, Hoboken, NJ 07030},

\cortext[cor1]{Corresponding author, Email address: alaa.ahmed@cornell.edu}

\begin{abstract}
\noindent Among the different concepts for wave energy conversion, oscillating surge wave energy converters have been shown to have a high capture width ratio. The primary wave capture structure consists of a flap hinged at the seabed or to a floating platform. Different flap configurations, including single and dual-flap, have been investigated. The separation distance between the oscillating surge wave energy converters can have an impact on their response when deployed in arrays. We consider the case of a dual-flap oscillating surge wave energy converter and investigate the impact of the separation distance between them on the performance of each flap. We estimate the absorbed wave energy and the annual energy production by the two flaps when deployed at the PacWave South site. Inviscid numerical simulations were conducted to predict the response of the oscillating surge wave energy converters. The simulations are validated with experimental measurements of a 1:10 scaled model in a wave tank. The results show that for a short separation distance, the interaction between the oscillating surge wave energy converters has a destructive and constructive effect depending on the wave frequency. However, these effects tend to balance each other out when considering the broad range of wave excitations. For longer separation distances, the interaction always results in a constructive effect. The results reveal that the separation distance has an insignificant impact on annual energy production when considering all wave frequencies and amplitudes.
\end{abstract}

\begin{keywords}
\texttt{Wave energy \sep Wave tank tests \sep OSWEC \sep WEC array \sep Separation distance \sep Power matrix}
\end{keywords}

\maketitle

\section{Introduction}

Wave energy converters (WEC) are often arranged in arrays to generate substantial amount of energy and serve as an effective alternate source for energy. When a WEC is deployed in water, it interacts with the incident waves and becomes a source for diffracted and radiated waves. Both diffracted and radiated waves have the characteristics of carrying energy similar to the incident waves. However, they can propagate in all directions from the WEC. The resultant fluid-structure interaction in a WEC array could affect all the other WECs in all directions. Consequently, this impacts the performance of each WEC in a different way depending on its location in the array, distance from other WECs, geometry, and wavelength of incident waves. The impact can be constructive in a way it adds to the energy absorption or destructive in that it decreases it \cite{budal1977theory}. Thus, assessing the impact of the separation distance is important in designing a WEC array. One of the most efficient WEC designs is the oscillating wave energy converter (OSWEC) \cite{babarit2015database,babarit2012numerical}. It consists of a single degree of freedom flap, hinged from one end at the bottom and allowed to rotate at the other end as the waves pass by. The relative motion between the flap and the waves is converted to power by using a power take-off (PTO) system \cite{li2021characterization}. It can be hinged directly to the seabed in shallow water or to a floating platform in deep water, which in turn is connected to the seabed by mooring lines \cite{mi2022dual}. \citet{folley2007effect} noted that positioning the OSWEC in deeper waters will result in larger captured power. Therefore, building an array of OSWECs in deep water can generate a substantial amount of energy. To be able to understand the behavior and response of the OSWEC in an array, the dynamic response should be examined with different separation distances. Reducing the distance between the flaps could result in a significant reduction in the mass of the platform and, consequently, the levelized cost of electricity (LCOE). However, it might add to the mooring loads and require a modified mooring configuration \cite{yang2020wave}. Therefore, other factors should be considered when designing an array and setting the distance between the OSWECs, such as mooring loads and manufacturing costs, although this is not the scope of the presented study.

\citet{budal1977theory} was the pioneer in the research field of studying the wave interaction between WECs in arrays. He demonstrated analytically the possibility of raising the power absorbed from identical bodies to 100\% of the incident wave. It can be accomplished when the bodies are positioned with a separation distance of a certain ratio to the wavelength. Several studies adopted the same approach with different designs but equally spaced in the array. The results highlighted the occurrence of constructive interference when a suitable layout is sought \cite{falnes1980radiation,thomas1981arrays,maniar1997wave}. Further experimental studies were performed under regular and irregular wave conditions to maximize the captured power by a single WEC in a small array (5 to 10 WECs). The experiments demonstrated that maximizing the power depends on the incident wave period and number of WECs in the array \cite{weller2010experimental,sun2016linear,stansby2017large}. Exploring this interaction experimentally with WECs in full-scale can be challenging and expensive; therefore, numerical tools are utilized to predict the hydrodynamic interaction between the WECs. Extensive research using the boundary element method (BEM) has been performed to numerically solve the hydrodynamic response of the WECs with consideration of the wave interference for arbitrary designs, the number of WECs in the array, and different ranges for separation distance. \citet{babarit2010impact} studied the effect of long separation distances on the energy production of a farm, which consists of two WECs, heaving cylinders and surging barges. He found the impact of the wave interference still affects the performance with a separation distance of 20 km for the surging barges. In both cases, a constructive effect was found on energy production. Moreover, the study considered the impact of the wave direction on the total energy production of the array. It was found that a zero heading angle maximizes energy production, while a 90-degree heading angle does not allow for any energy production. \citet{borgarino2012impact} considered the same WECs as in \cite{babarit2010impact}, but for larger arrays of 9–25 WECs and different layouts for the farm. The research concluded that the rear rows in the farm will experience destructive effects. However, the effect on mean energy production is minor. Further numerical and analytical studies on heaving and surging WECs concerning the wave interference and its effect on the performance of the WECs can be found in \cite{falcao2002wave,ricci2007point, child2007interaction,vicente2009dynamics,child2009modification, babarit2010assessment,de2010power,child2010optimal,goteman2015optimizing, mcguinness2016hydrodynamic,goteman2017wave,2017,moarefdoost2017layouts, kara2019maximise,stansby2022total,goteman2022multiple,liu2022wave}.

Most of the aforementioned research is related to the wave interference within arrays of WECs with transitional displacement, while the interest in this paper is focused on OSWECs and their performance in an array. \citet{babarit2013park} studied the park effect in arrays of OSWECs and found that for small arrays with conventional sizes (dimensions of 10 to 20m) with common layouts with separation distances of 100 to 200m, the parking effect will be negligible. On the other hand, for large arrays (more than 10 OSWECs), a negative park effect is increasing with the increasing number of rows. His findings recommend laying out the array with a small number of rows as much as possible, which come in agreement with the findings in \cite{noad2015optimisation}. Furthermore, \citet{renzi2014wave} investigated the dynamics of a wave farm consisting of two OSWECs in the nearshore. They identified a near-resonant phenomenon, which is a distinctive property of an array of OSWECs, resulting in constructive interference in the array. Although the OSWEC is a promising design for WECs with higher efficiency in capturing power, the analysis of the finite array of OSWECs has been limited and hasn't been explored broadly yet.

To date, all of the research conducted on wave interference within arrays of WECs is based on the assumption of potential flow and linear wave theory with small wave amplitudes. In real sea-states, these assumptions are no longer valid, and a more accurate approach is needed to include all the nonlinearities in the wave field surrounding the WECs. High-fidelity simulations solving Reynolds-averaged Navier-Stokes (RANS) equations can give accurate insights, but they are computationally expensive. Alternatively, medium-fidelity simulations solving Euler's equations can balance between computational time and accuracy. Throughout this paper, a full-scale dual-flap OSWEC \cite{AHMED2024130431} is considered to serve as an array of two OSWECs. Utilizing medium-fidelity simulations, the impact of the separation distance between the flaps on the performance of the flaps and the annual energy production is assessed under real sea-state conditions. Torque-forced simulations were utilized to assess the impact of wave interaction, radiation, and diffraction on the response of each flap. Wave-forced numerical simulations were conducted to investigate the effect of the location of each flap relative to the incident waves and the total impact of each flap on the other's performance. Variant wave heading angles were tested to study the effect of wave direction on the performance of the flaps. Finally, a power matrix is estimated for each tested separation distance to evaluate its impact on the annual energy production.

\section{Considered WEC design }

\begin{figure}
         \centering
         \includegraphics[clip,width=0.9\textwidth]{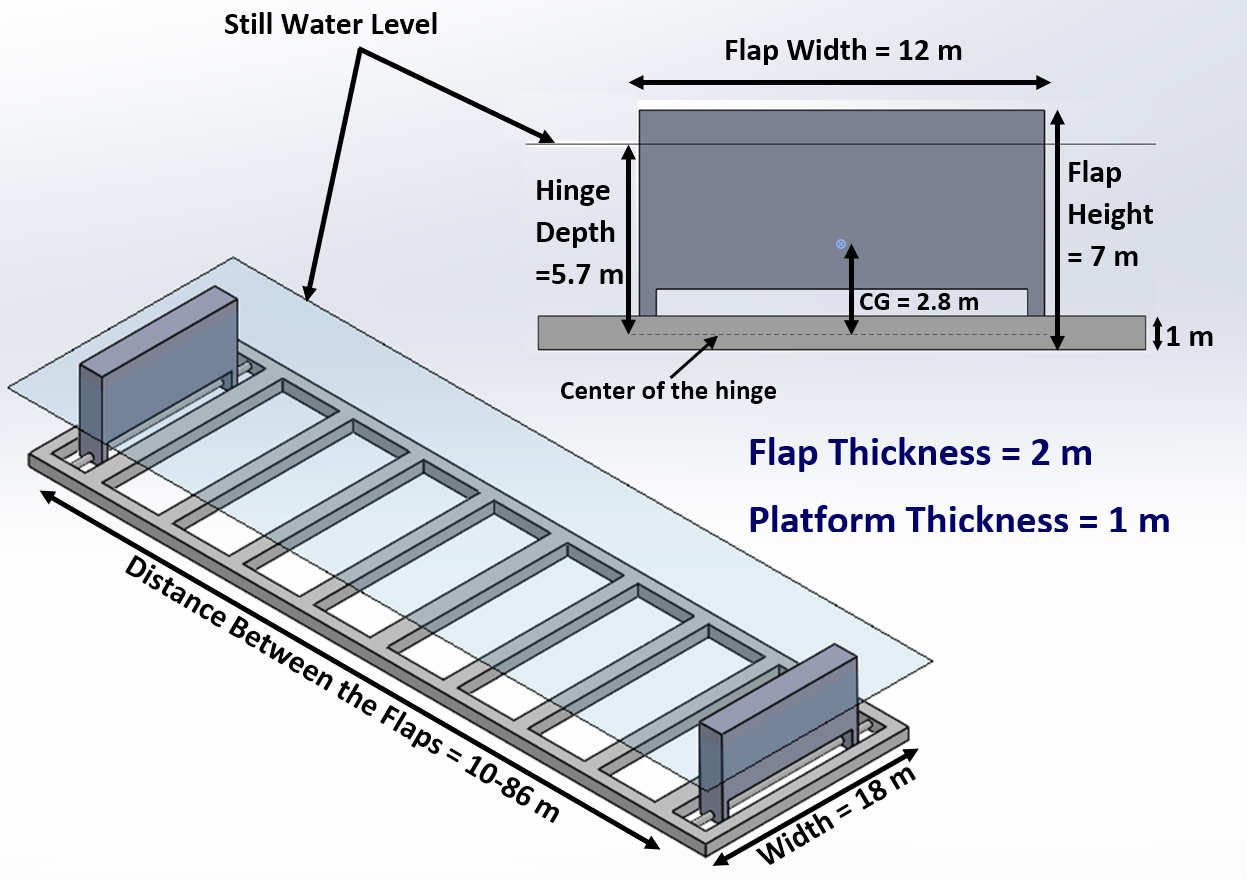}
         \caption{Schematic of full-scale dual-flap OSWEC}
         \label{fig: full-scale}
\end{figure}

We considered the same design of a dual-flap OSWEC as in \cite{AHMED2024130431} to study the effect of the separation distance between the flaps and optimize it to maximize the annual energy production and reduce mooring loads and manufacturing costs. The dual-flap OSWEC, as schematically shown in Figure \ref{fig: full-scale}, can serve as an array of two devices. However, small separation distances ranging from 10 m to 86 m are considered in this study. Each flap is hinged to a horizontal shaft and positioned perpendicular to the incident wave. In this paper, the platform is kept fixed, and the system is assessed as being considered to be deployed at the PacWave South site \cite{dunkle2020pacwave}.

A 1:10 model of the proposed dual-flap OSWEC is fabricated and tested in the wave tank of the Davidson Laboratory at Stevens Institute of Technology. The experimental data was used for numerical validation. More details about the experimental setup and numerical model validation can be found in \cite{AHMED2024130431,AHMED}

\section{Methodology}

A WEC in an array can be a source of diffracted and radiated waves as a result of its interaction with incident waves. Both diffracted and radiated waves can impact the performance of the other WECs in the array. To understand the wave interaction between the WECs in an array, the dynamic response of each WEC should be studied. The interaction and impact of each flap on the others can be split into two main reasons: how the front WEC is interacting with the incident waves before reaching the back WECs, and how its response is contributing to the waves reaching the back ones. To isolate each effect and be able to visualize the impact of each one independently, we tested the considered dual-flap OSWEC under torque-forced and wave-forced. The torque forcing can allow us to isolate each forcing component contributing to the motion to be able to investigate its effect independently. While the wave forcing will show the impact of the front flap on the surrounding wave field and how the field is changing before reaching the back flap, it will also give a realistic estimate for the annual energy production at the targeted site. Experimental testing for this study is expensive and challenging; therefore, numerical simulations can play a significant role. Considering the number of simulations required for a comprehensive investigation, RANS simulations will be inefficient. The same approach in \cite{AHMED2024130431} was followed by utilizing inviscid simulations.

\subsection{Torque-forced simulations}
Torque-forced simulations can help in isolating the impact of each component in the equation of motion of the system and investigating its impact on the response of each OSWEC. The response of a single flap is governed by applying Newton's second law, yielding the equation of motion to:
\begin{equation}
(I+I_a)\ddot \theta + C \dot \theta + k \theta = T_{_0} sin(\omega t + \phi)
\label{eq: eq motion der}
\end{equation}

while for the dual-flap OSWEC, the coupling and interaction between the responses of the flaps should be considered. The equation of motion then yields:
\begin{equation}
\begin{split}
\left(
\begin{bmatrix}
I_{ll} & 0\\
0 & I_{rr} 
\end{bmatrix}
+
\begin{bmatrix}
I_{a_{ll}} & I_{a_{lr}}\\
I_{a_{rl}} & I_{a_{rr}} 
\end{bmatrix}
\right)
\begin{bmatrix}
\ddot \theta_l\\
\ddot \theta_r 
\end{bmatrix}
+
\begin{bmatrix}
C_{ll} & C_{lr}\\
C_{rl} & C_{rr} 
\end{bmatrix}
\begin{bmatrix}
\dot \theta_l\\
\dot \theta_r 
\end{bmatrix}
\\+
\begin{bmatrix}
k_l\\
k_r 
\end{bmatrix}
\begin{bmatrix}
\theta_l\\
\theta_r 
\end{bmatrix}
=
\begin{bmatrix}
T_{_{0l}} sin(\omega t+\phi_l)\\
T_{_{0r}} sin(\omega t+\phi_r)\\
\end{bmatrix}
\end{split}
 \label{eq: EoM}
\end{equation}

where $I$ is the mass moment of inertia of the flap, $I_a$ is the added mass moment of inertia, $C$ is the damping coefficient, and $k$ is the stiffness coefficient. $\theta$, $\dot\theta$, and $\ddot \theta$ are the angular displacement, velocity, and acceleration, respectively. $I_{a_{lr}}$ and $I_{a_{rl}}$ are the coupling mass moment of inertia between the flaps; $C_{lr}$ and $C_{rl}$ are the coupling damping coefficients; $\phi$ is the phase difference angle; and $T_{_0}$, $T_{_{0l}}$, and $T_{_{0r}}$ are the excitation torque amplitudes.

It should be noted that $I=I_{ll}=I_{rr}$, $I_a=I_{a_{ll}}=I_{a_{rr}}$, $I_{a_{lr}}=I_{a_{rl}}$, $C=C_{ll}=C_{rr}$, $C_{lr}=C_{rl}$, and $k=k_l=k_r$ for the same separation distance and excitation period.

The impact of changing the separation distance between the flaps is investigated under harmonic torque-forced with different excitation periods and amplitudes. Three separation distances ($d=10\: m, 45\: m, \:and\: 70\: m$) and four excitation periods ($T_e = 7.5 \: s, 8.5 \: s, 9.5 \: s, \:and \: 10.5 \: s$) are selected for this case study based on the targeted deployment site.

Torque-forced simulations are conducted by feeding the excitation torque as an input to each flap and observing its dynamic response. To investigate the effect of each term in the equation of motion, we tested different excitation cases where the amplitude of the excitation torque and the phase differences were changed. To cancel the effect of the location of each flap in the flow field, which might cause additional damping or excitation as in the case with incident waves, the input excitation torque is applied with the same amplitude on both flaps. The flaps are noted by the left flap and the right flap as the location effect is eliminated from the response. In the first case, the excitation torque is applied only to the right flap, and the left flap is kept fixed. This case can show the impact of the distance on the response of the right flaps due to the reflected and diffracted waves from the left fixed flap, which is implicitly represented in the excitation term in Equation \ref{eq: EoM}. In the next case, the excitation torque is still applied to the right flap only, but the left flap is allowed to rotate as a response to its interaction with the perturbations in the wave field caused by the radiated waves generated from the right one. This case assesses the impact of the coupling terms in the equation of motion on the response of the flaps. Finally, three other cases are considered for both flaps being excited, but with varying phase differences ($\phi$) to assess the impact of the phase difference on the interaction between the flaps. Three phase differences are tested; the first two are exciting the flaps in-phase and out-of-phase, which correspond to $\phi=0^0 \:and \:180^0$, respectively. The last phase difference, namely, the arbitrary phase, represents the ratio between the separation distance between the flaps ($d$) and the wavelength ($\lambda$) of the excitation period ($\phi=\frac{d}{\lambda}2\pi$). The last case represents a more realistic one, as it would be with the deployment of the flaps in the ocean. All five cases are compared to a single flap case to observe the change in response.

\subsection{Wave-forced simulations}
The impact of changing the separation distance between the flaps is investigated under regular wave forcing with different excitation periods and amplitudes. In order to assess its impact on annual energy production, a power matrix can be estimated for each distance for the wave resources at the PacWave south site and compared all together. To save computational power, a partial power matrix was estimated using the highest frequency-occurring waves and the waves with the highest energy percentage. These waves have periods between 7.5 s and 11.5 s and heights between 1.25 m and 5.25 m. Numerical simulations solving Euler equations were conducted. Seven different separation distances were tested to cover a range of 10 m to 86 m. In one case, the distance was set to half the wavelength of the wave with the highest percentage of total energy (86 m) based on the wave resources at the targeted site. In another case, we tested the performance when the distance was set to half the wavelength of the wave having an equal frequency to the natural frequency of the flap (70 m). The third distance was set to half the wavelength of the most-occurring wave (55 m). Finally, we reduced the distance to half the wavelengths of four short waves (45 m, 33 m, 15 m, and 10 m) to reduce the length and mass of the floating platform. Setting the distance to be half the wavelength of the occurring waves can exploit the out-of-phase motion between the flaps and reduce the mooring loads.

\section{Results}
Torque-forced and wave-forced simulations were validated against experimental data from wave tank tests on a 1:10 scaled-model of the OSWEC. The validation of numerical simulations is found in \cite{AHMED2024130431,AHMED}. The results show that FLUENT can predict the response of the OSWEC by solving Euler’s equations under different conditions with acceptable accuracy, which concludes the reliability of using inviscid numerical simulations to investigate the impact of the separation distance between the flaps on the annual energy production. The same approach utilizing inviscid simulations is adopted for all the simulations presented in this paper on the full-scale dual-flap OSWEC.

\subsection{Torque-forced cases}
\begin{figure}
\centering
\subcaptionbox{}{\includegraphics[trim={0.5cm 0.1cm 1cm 0cm},width=0.45\textwidth]{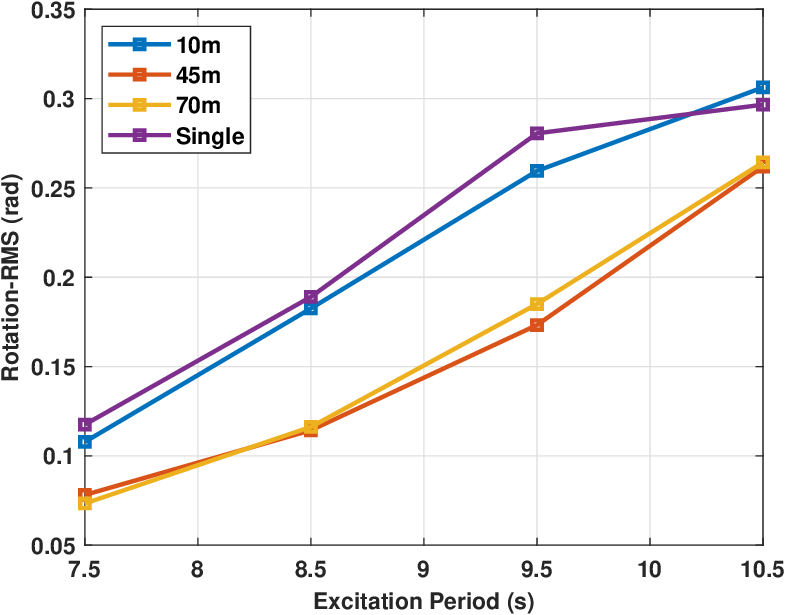}}%
\hfill
\subcaptionbox{}{\includegraphics[trim={0.5cm 0.1cm 1cm 0cm},width=0.45\textwidth]{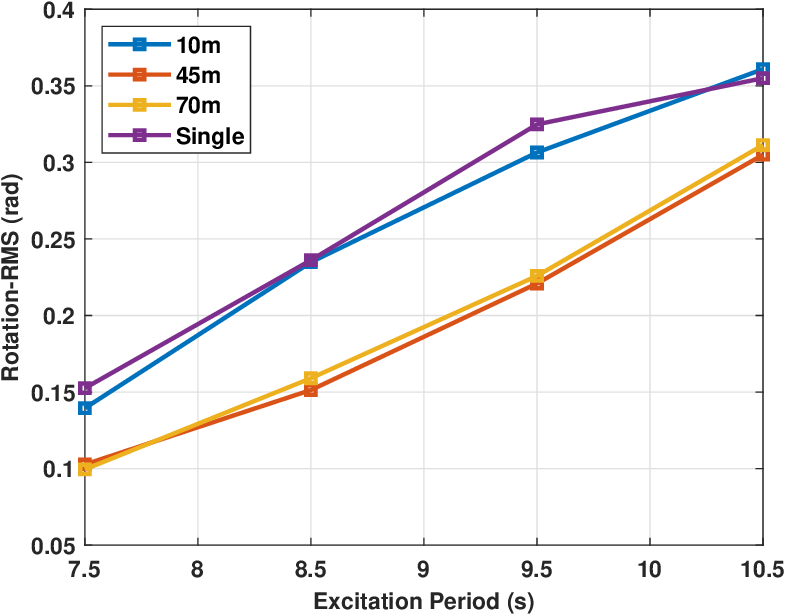}}%
\vfill
\subcaptionbox{}{\includegraphics[trim={0.5cm 0.1cm 1cm 0cm},width=0.45\textwidth]{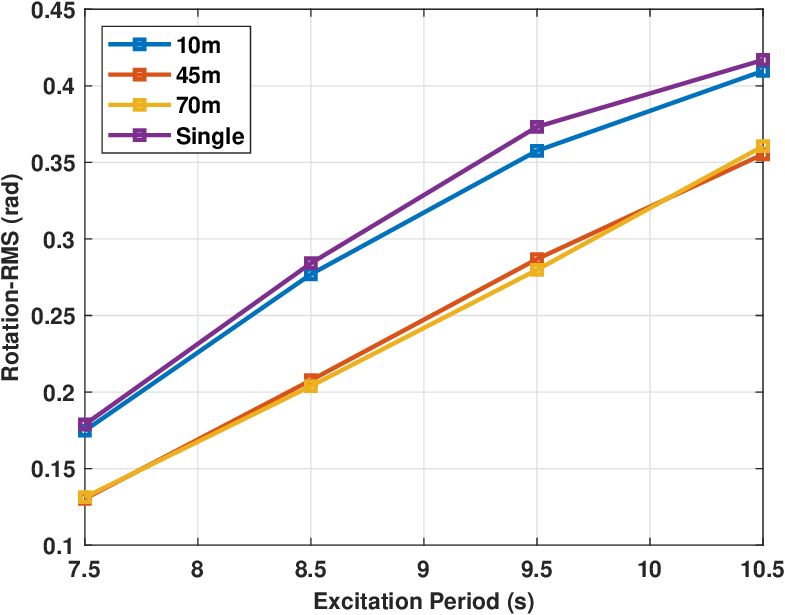}}%
\hfill
\subcaptionbox{}{\includegraphics[trim={0.5cm 0.1cm 1cm 0cm},width=0.45\textwidth]{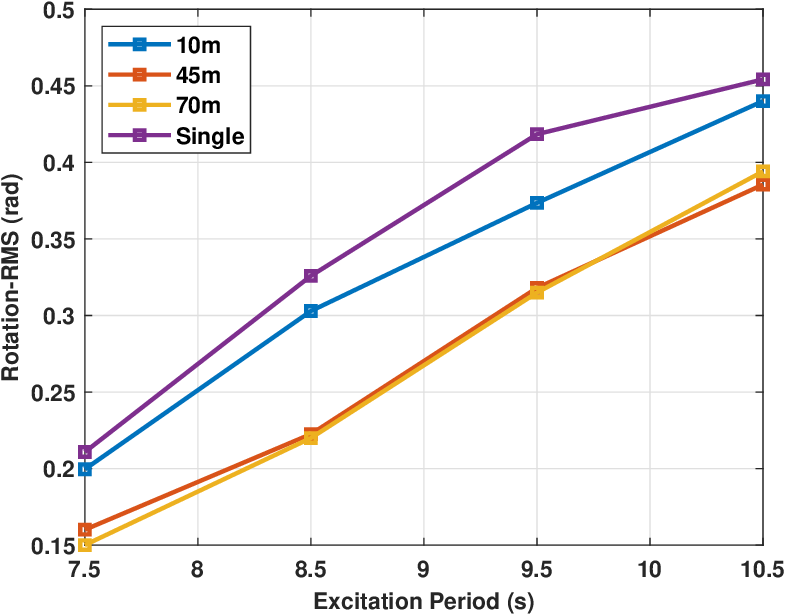}}%
\caption{The RMS values of the rotation of the right flap under regular forcing with excitation torque amplitudes of (a) 0.6 MN.m and (b) 0.8 MN.m (c) 1 MN.m (d) 1.2 MN.m with the left flap fixed}
\label{fig: fixed case}
\end{figure}

The right flap was excited by different amplitudes for the excitation torque, ranging from 0.6 MN.m to 1.2 MN.m in the first two cases. Figure \ref{fig: fixed case} compares the RMS values of the rotation of the right flap to the single flap when the left flap is kept fixed. The results show that for all excitation periods and amplitudes with a separation distance of 10 m, the difference between the response of a single flap and the response of the right flap is minimal and negligible. This could be due to having equal constructive and destructive effects from radiated and diffracted waves that cancel the impact of the left flap. While for the distances of 45m and 70m, the destructive effect from diffraction is higher than the constructive effect from radiation, resulting in a reduction in the RMS value of the rotation compared to the single flap.  However, changing the distance from 45m to 70m has no effect on the response of the right flap. This concludes that when the distance is in the range of $0.06\lambda < d < 0.11\lambda$, the response of the flap is not affected by the presence of another one due to the compensation between the destructive and constructive effects. While for a longer distance in the range of $0.25\lambda < d < 0.80\lambda$, the destructive effect from diffraction becomes dominant on the response.

\begin{figure}[H]
\centering
\subcaptionbox{}{\includegraphics[width=0.35\textwidth]{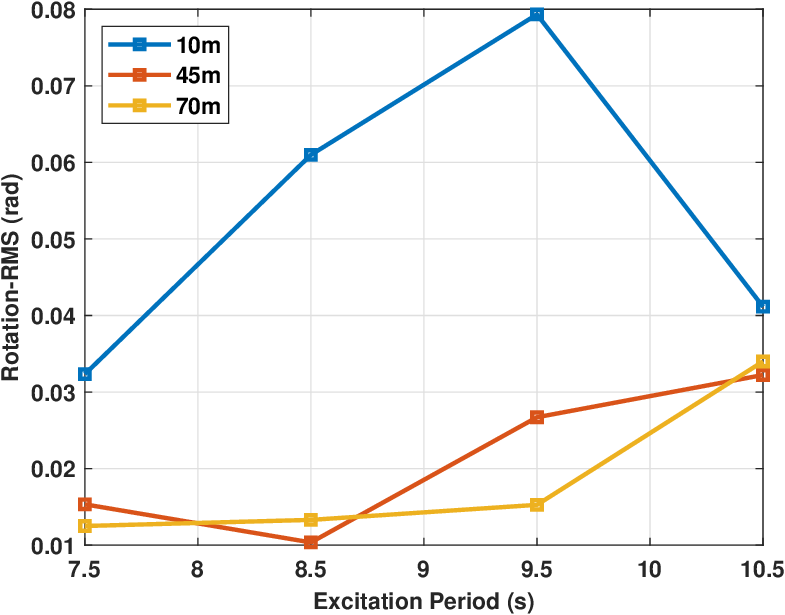}
\hfill
\includegraphics[width=0.35\textwidth]{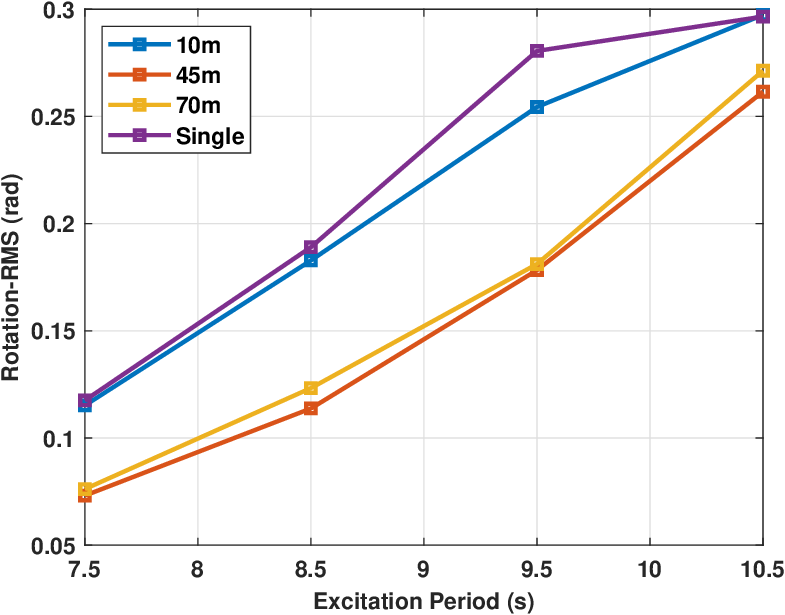}}%
\vfill
\subcaptionbox{}{\includegraphics[width=0.35\textwidth]{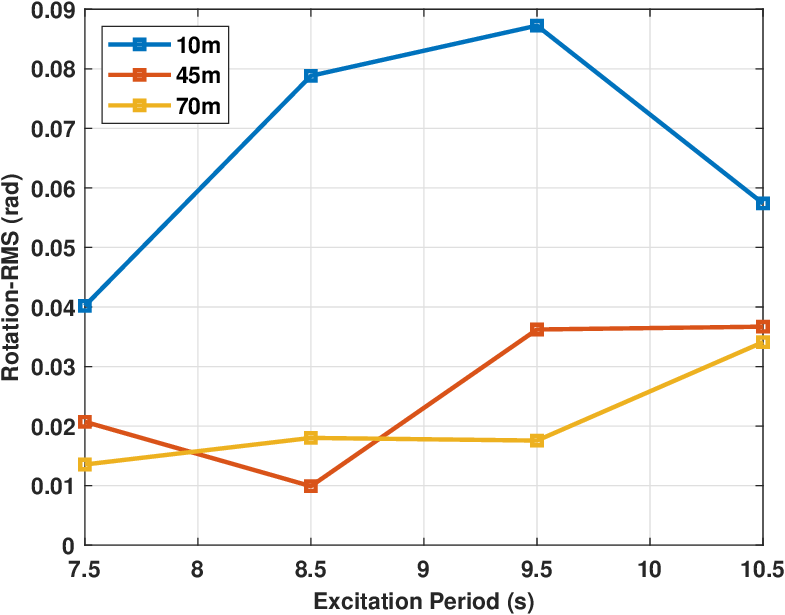}
\hfill
\includegraphics[width=0.35\textwidth]{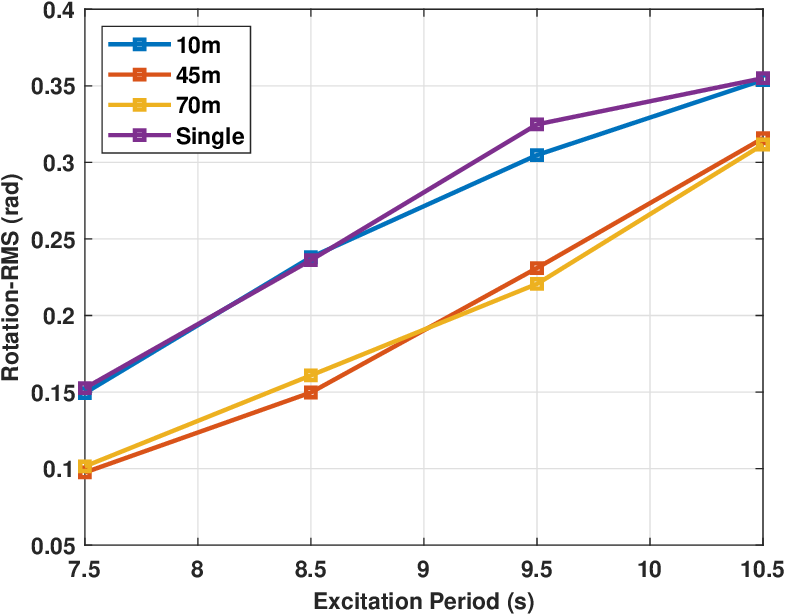}}%
\vfill
\subcaptionbox{}{\includegraphics[width=0.35\textwidth]{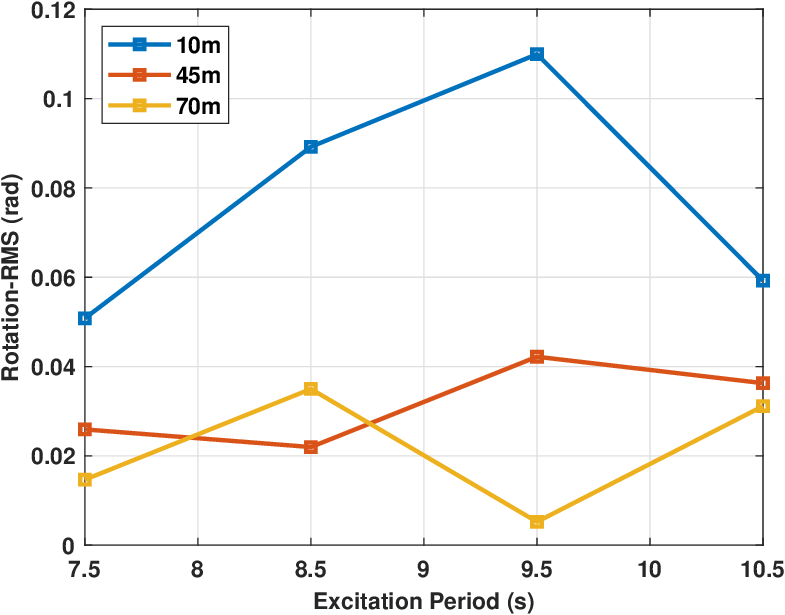}
\hfill
\includegraphics[width=0.35\textwidth]{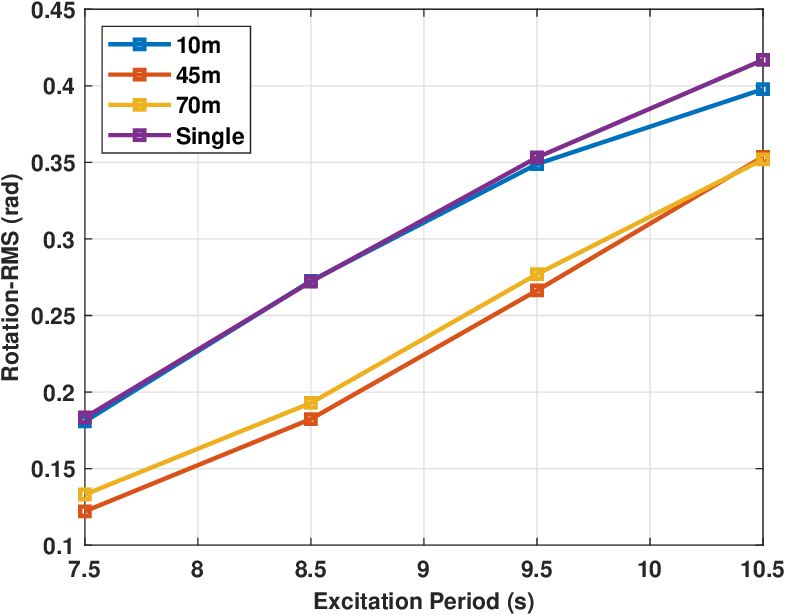}}%
\vfill

\subcaptionbox{}{\includegraphics[width=0.35\textwidth]{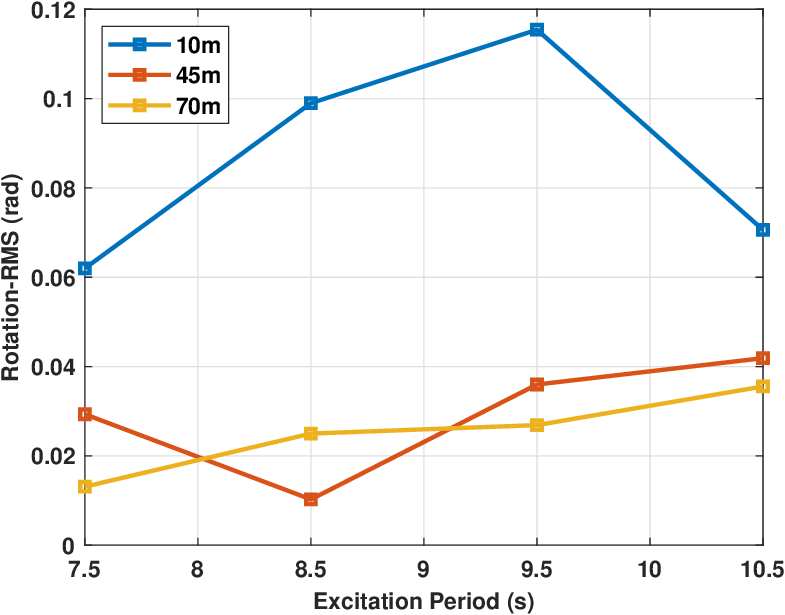}
\hfill
\includegraphics[width=0.35\textwidth]{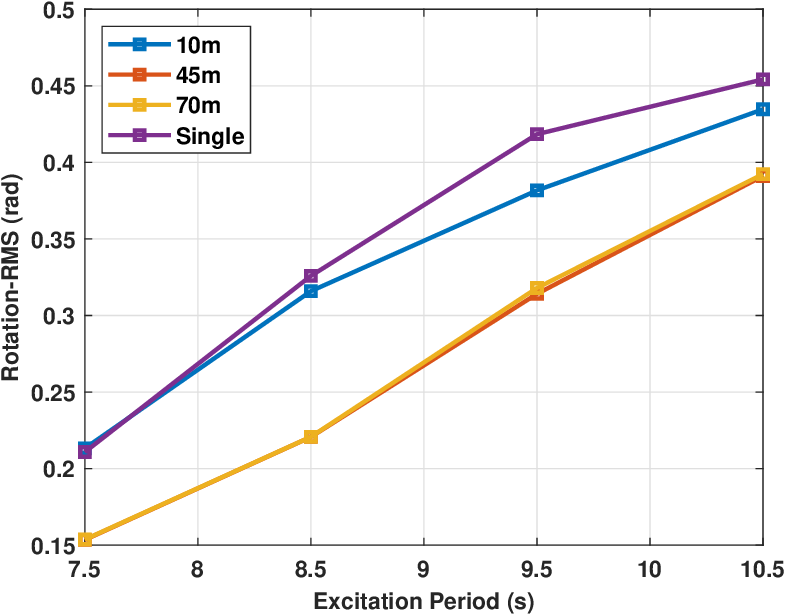}}%
\caption{The RMS values of the rotation of the left flap (left) and right flap (right) with single flap under regular forcing with excitation torque amplitude of (a) 0.6 MN.m (b) 0.8 MN.m (c)1 MN.m (d) 1.2 MN.m applied on the right flap}
\label{fig: rightexc}
\end{figure}

Figure \ref{fig: rightexc} shows the RMS values of the rotation of both flaps when the right flap is excited and the left flap is allowed to rotate as a response to its interaction with the wave field generated by the right flap. The right plots compare the rotation RMS of the right flap with the rotation RMS of the single flap. It can be shown that for all different excitation periods and all different amplitudes, the rotation RMS of the left flap increases as the separation distance decreases. This was expected as the wave field gets stronger with smaller distances, which strengthens the interaction between the flap and the surrounding wave field, generating a constructive effect. However, the rotation of the left flap is very small relative to the right flap and has no additional effect on the response of the right flap in this case when compared to the fixed case. This confirms the findings from Figure \ref{fig: fixed case}, that radiation has a small impact at longer distances. Findings from Figures \ref{fig: fixed case} and \ref{fig: rightexc} come in agreement with the findings in \cite{babarit2010assessment}. Given that the distances of 45 m and 70 m will cause some cases to be in-phase and some others to be out-of-phase, that explains the slight change noticed in the RMS values of the rotation of the left flap. This will be demonstrated more in the next case.

\begin{figure}[H]
\centering
\subcaptionbox{}{\includegraphics[width=0.5\textwidth]{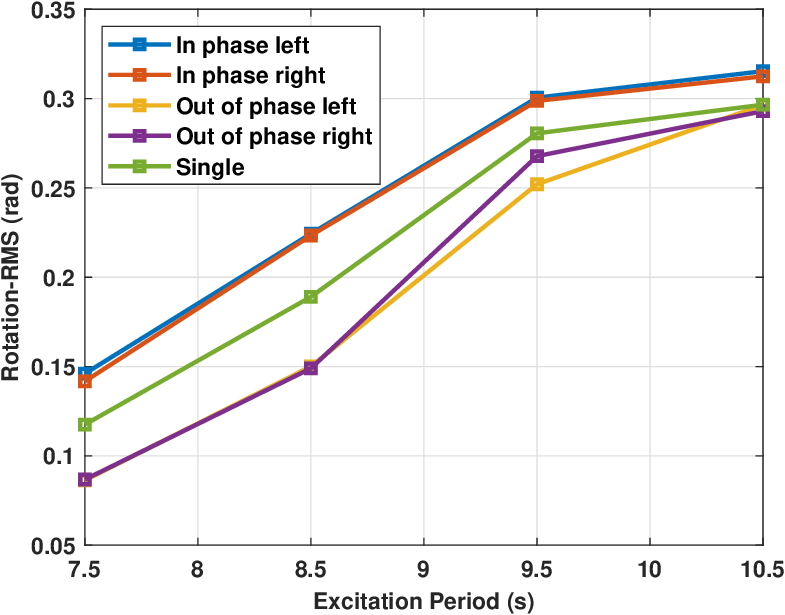}
\hfill
\includegraphics[width=0.5\textwidth]{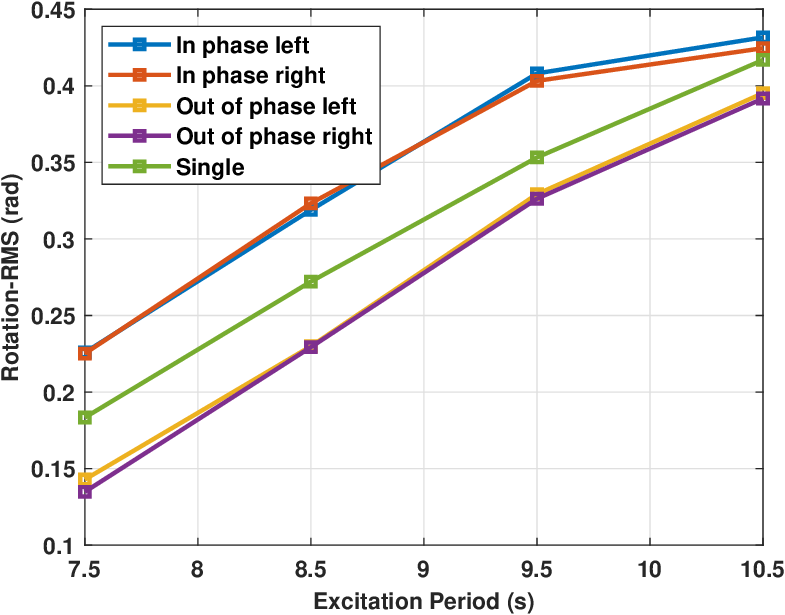}}%
\vfill
\subcaptionbox{}{\includegraphics[width=0.5\textwidth]{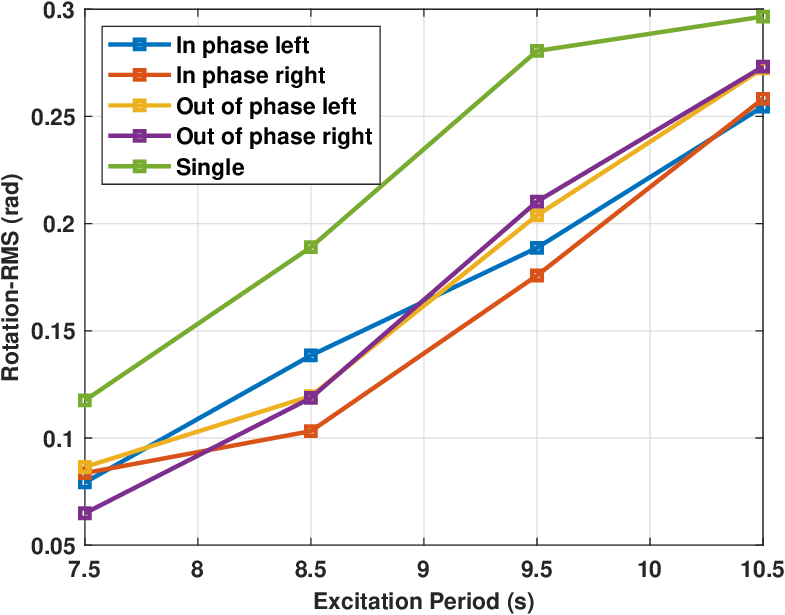}
\hfill
\includegraphics[width=0.5\textwidth]{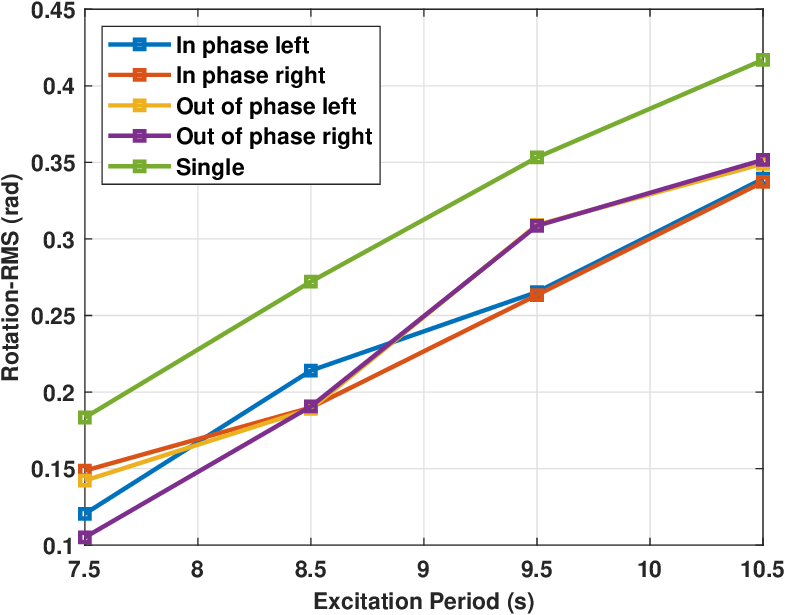}}%
\vfill
\subcaptionbox{}{\includegraphics[width=0.5\textwidth]{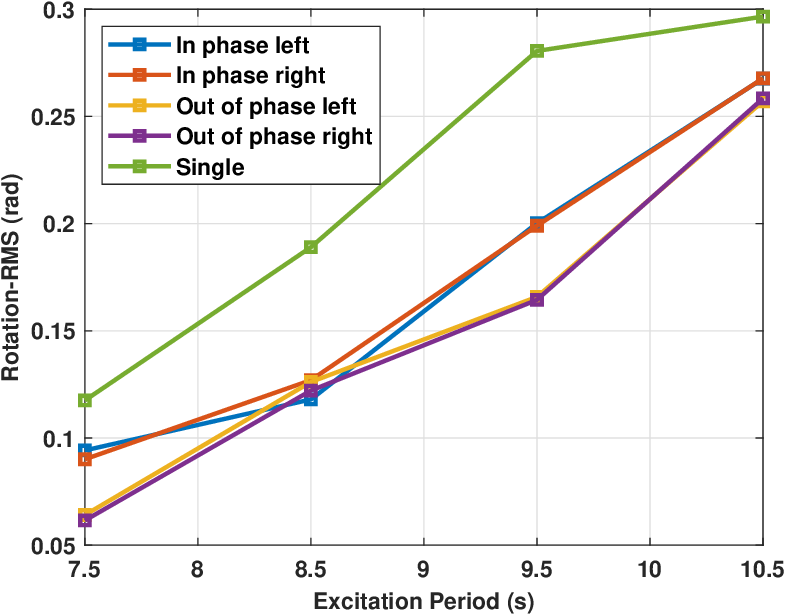}
\hfill
\includegraphics[width=0.5\textwidth]{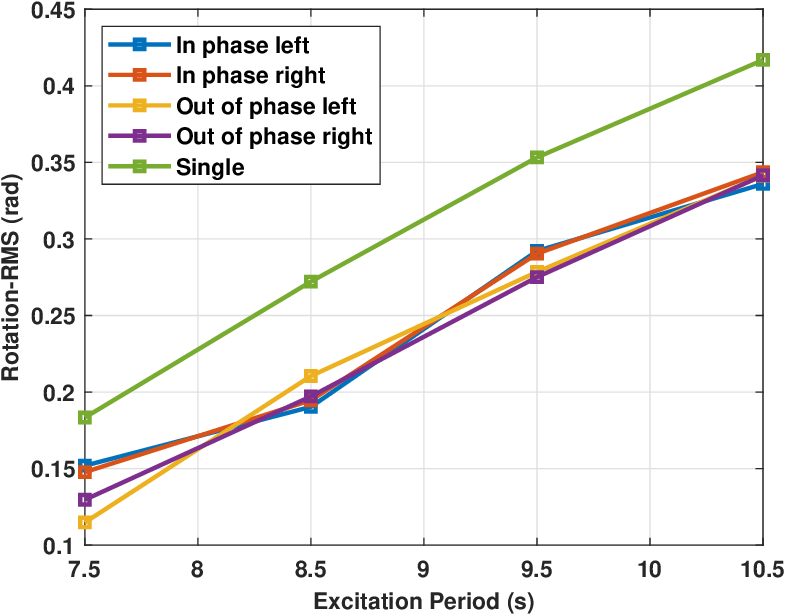}}%
\caption{The RMS values of the rotation of the flaps separated by (a) 10 m (b) 45m (c) 70 m under regular forcing with excitation torque amplitude of 0.6 MN.m (left) and 1.0 MN.m (right) amplitudes when the two flaps are in-phase and out-of-phase}
\label{fig: inout}
\end{figure}

Figures \ref{fig: fixed case} and \ref{fig: rightexc} show that changing the amplitude of the excitation torque doesn't change the trend of the responses of the flaps. Therefore, for further investigations, two amplitudes for the excitation torque will be considered: 0.6 MN.m and 1MN.m. To check the impact of the phase differences, both flaps were excited with varying phase differences between the flaps. The first configuration is when both flaps were excited in-phase ($\phi_l=\phi_r=0$); the second configuration is when they were excited out-of-phase ($\phi_l=0, \phi_r=180^0$). Figure \ref{fig: inout} shows that for a separation distance of 10 m, the RMS values of the rotation of both flaps are higher than the single flap for in-phase configurations, while lower in the out-of-phase configuration. Due to the small separation distance, the impacts of the radiated and diffracted waves on the flaps during the in-phase configuration are adding up, which results in additional excitation from one flap on the other, which contributes to a higher rotation. While in the out-of-phase configuration, the impacts are opposite, which results in additional damping from one flap on the other, yielding a lower rotation than the single flap.
For larger distances, the discrepancies in the responses between the two configurations are due to the resultant phase difference between the responses. However, these changes are insignificant over different excitation frequencies. The results show that for short distances with a range of $0.06\lambda < d < 0.11\lambda$, the phase difference has a significant impact on the response of the flaps. While for larger distances with a range of $0.25\lambda < d < 0.80\lambda$, the phase difference has minimal impact on the responses of the flaps.

\begin{figure}[H]
\centering
\subcaptionbox{}{\includegraphics[width=0.5\textwidth]{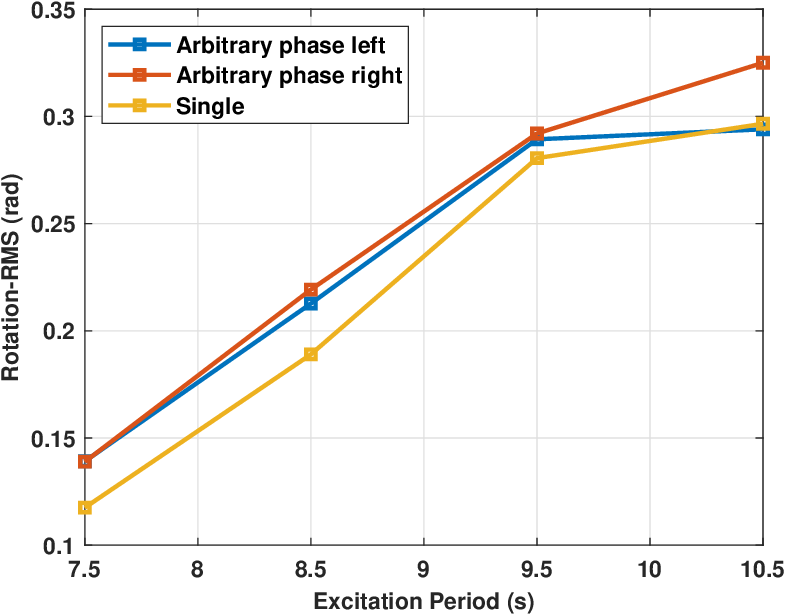}
\hfill
\includegraphics[width=0.5\textwidth]{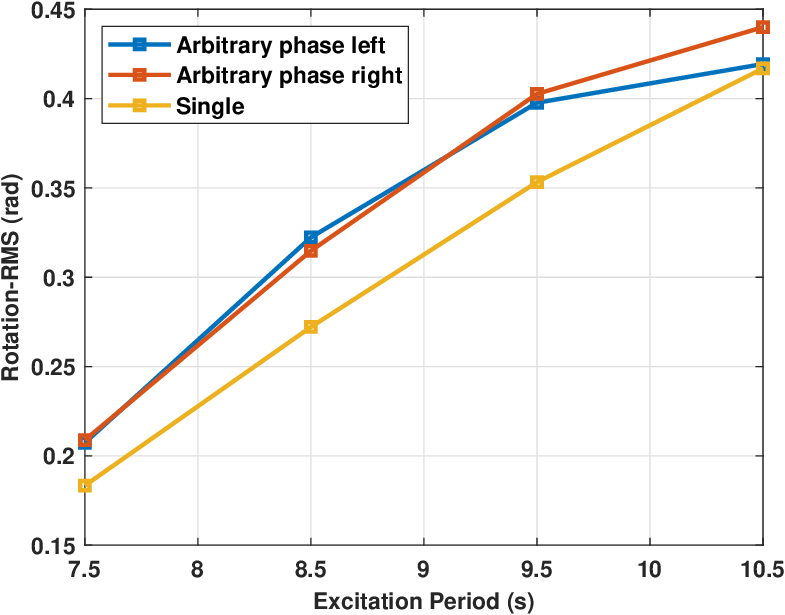}}%
\vfill
\subcaptionbox{}{\includegraphics[width=0.5\textwidth]{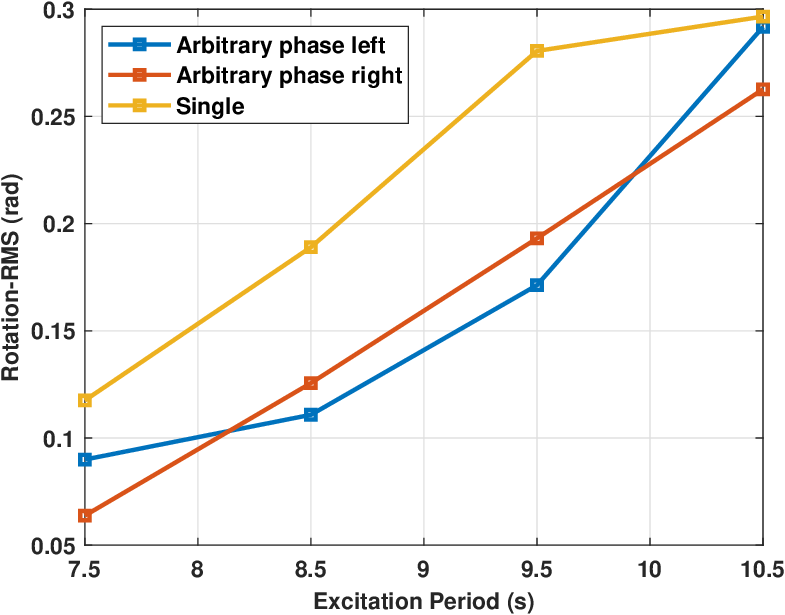}
\hfill
\includegraphics[width=0.5\textwidth]{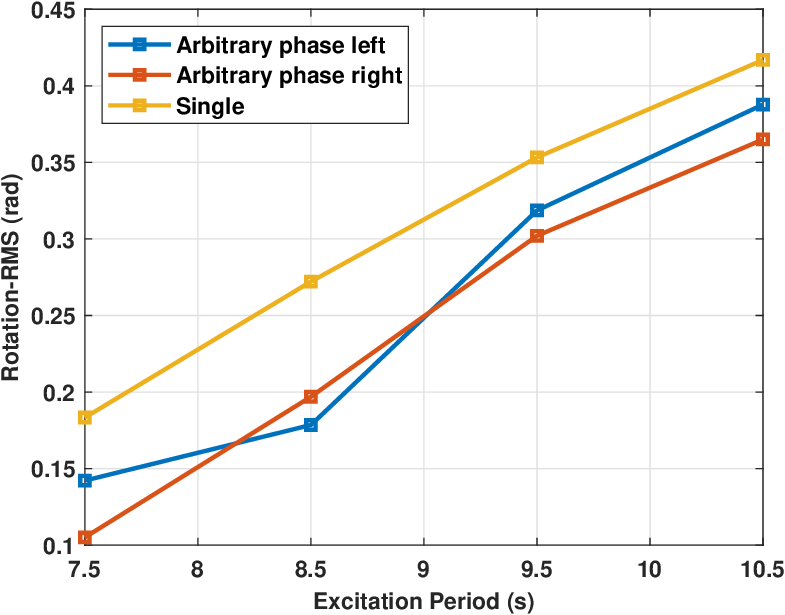}}%
\vfill
\subcaptionbox{}{\includegraphics[width=0.5\textwidth]{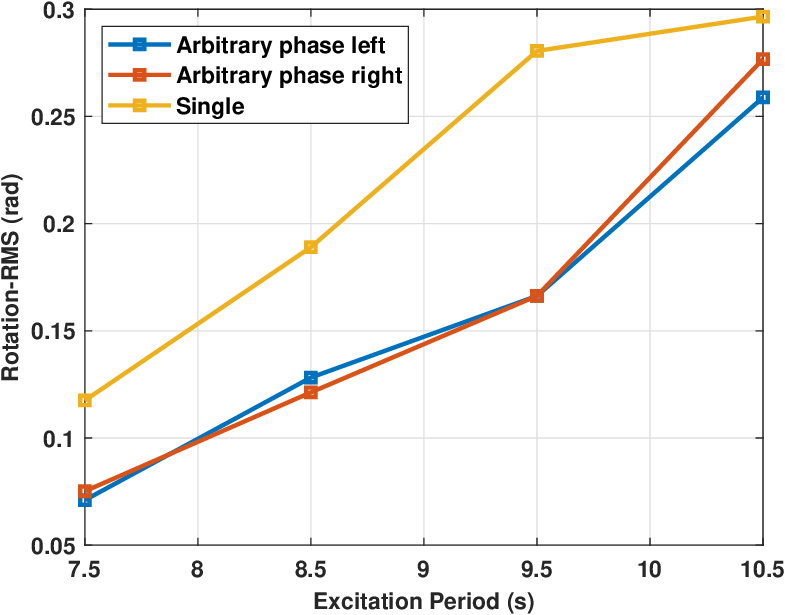}
\hfill
\includegraphics[width=0.5\textwidth]{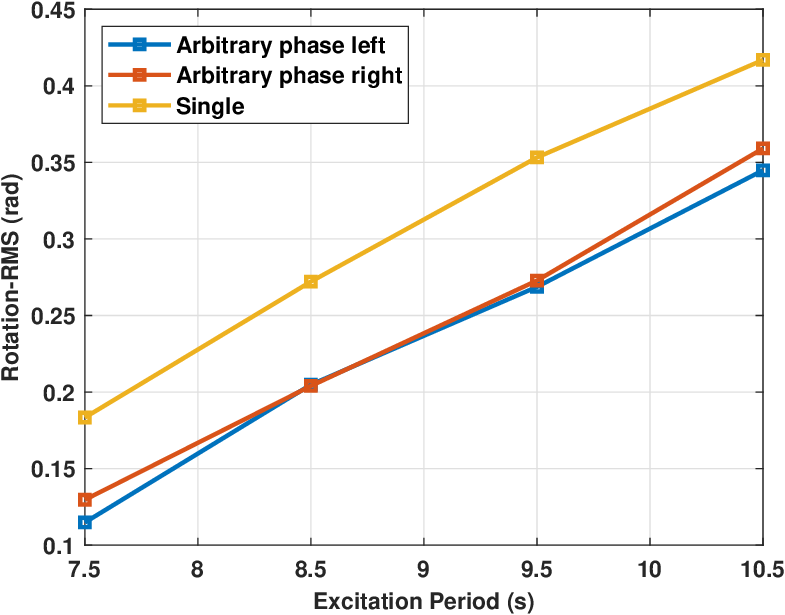}}%
\caption{The RMS values of the rotation of the flaps separated by (a) 10 m (b) 45m (c) 70 m under regular forcing with excitation torque amplitude of 0.6 MN.m (left) and 1.0 MN.m (right) amplitudes when the two flaps have arbitrary phase difference ($\phi=\frac{d}{\lambda}2\pi$)}
\label{fig: arbtr}
\end{figure}

To test a realistic excitation case as in the ocean and to be able to compare the results with a case of incident waves, the last configuration is presented with equal amplitudes of excitation torque applied on both flaps with an arbitrary phase ($\phi_l=0, \phi_r=\frac{d}{\lambda}2\pi$). Figure \ref{fig: arbtr} shows an increase in the RMS values of the rotation when the flaps are close to each other compared to the RMS values of the single flap, while the RMS values of the rotation are slightly decreased for larger distances. That confirms the conclusion of the significant impact of the phase differences on the response of the flaps with a short separation distance ($0.06\lambda < d < 0.11\lambda$) and the dominant impact of diffraction at longer distances ($0.25\lambda < d < 0.80\lambda$).

\FloatBarrier

\subsection{Wave-forced cases}
\begin{figure}
\centering
\subcaptionbox{}{\includegraphics[width=0.5\textwidth]{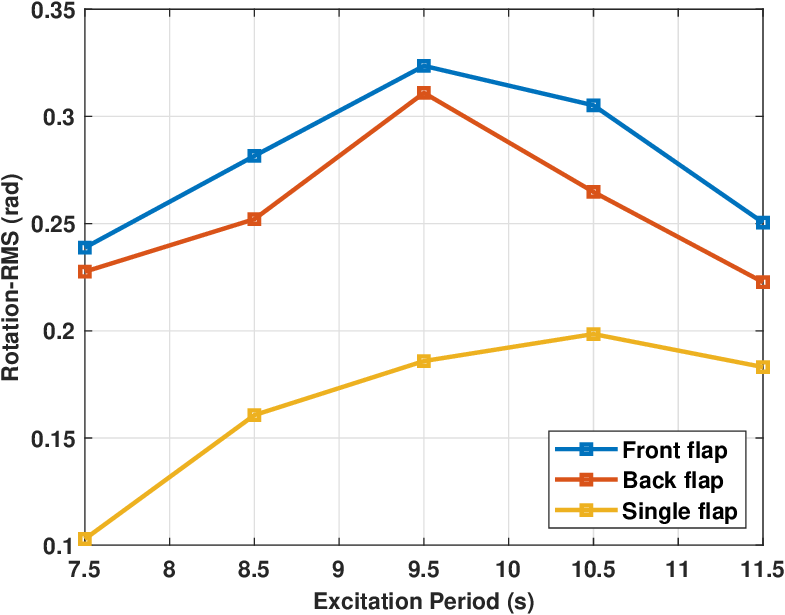}
\hfill
\includegraphics[width=0.5\textwidth]{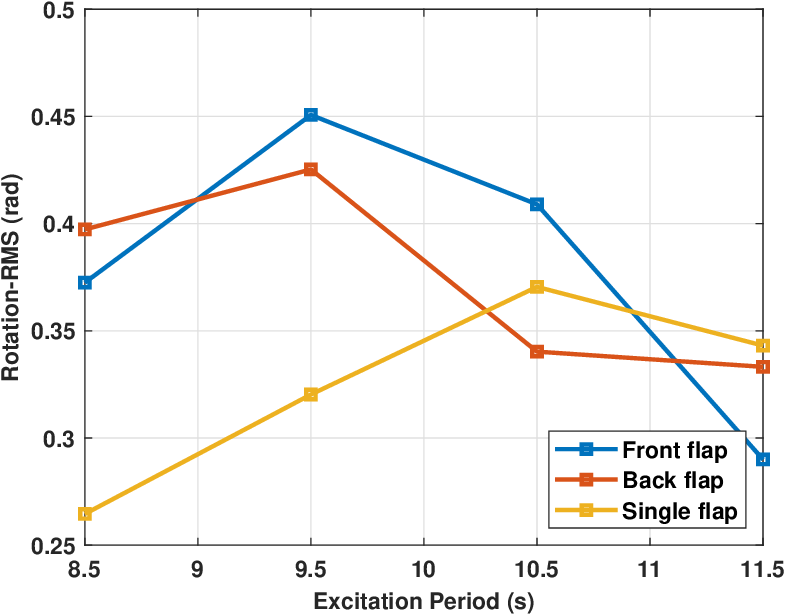}}%
\vfill
\subcaptionbox{}{\includegraphics[width=0.5\textwidth]{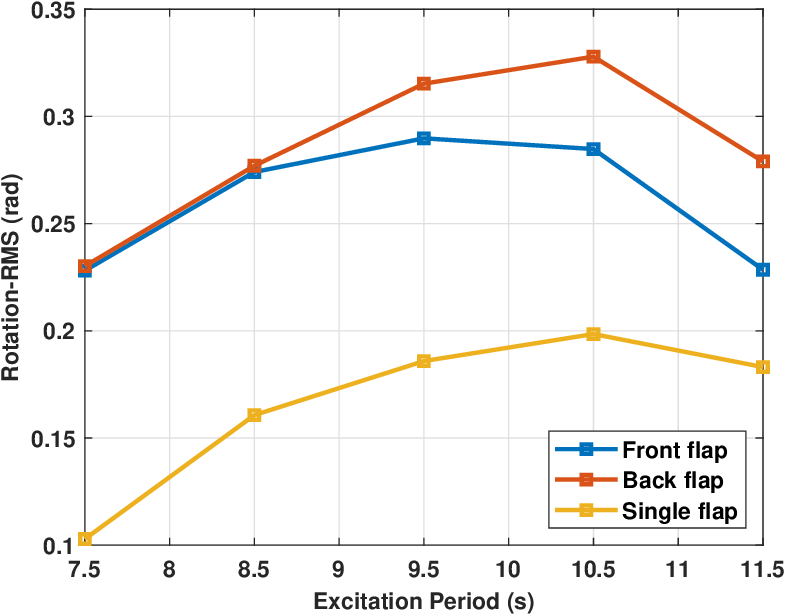}
\hfill
\includegraphics[width=0.5\textwidth]{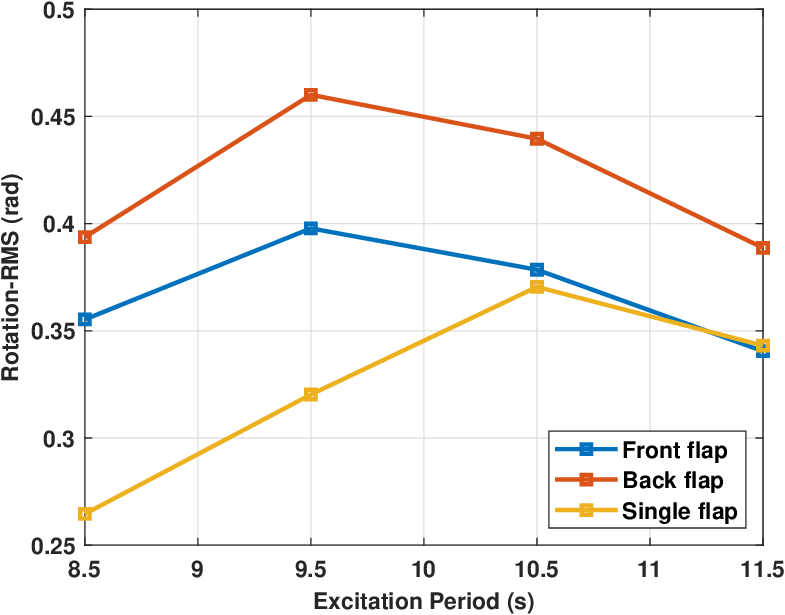}}%
\vfill
\subcaptionbox{}{\includegraphics[width=0.5\textwidth]{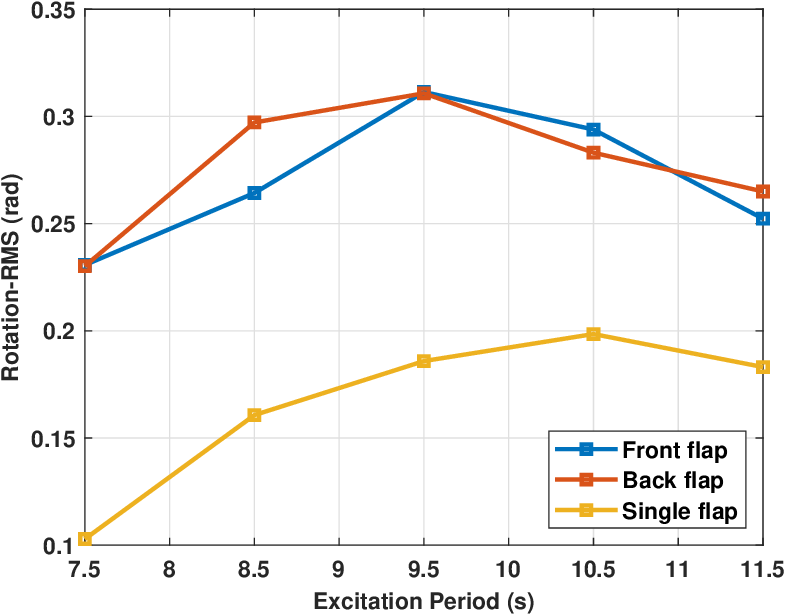}
\hfill
\includegraphics[width=0.5\textwidth]{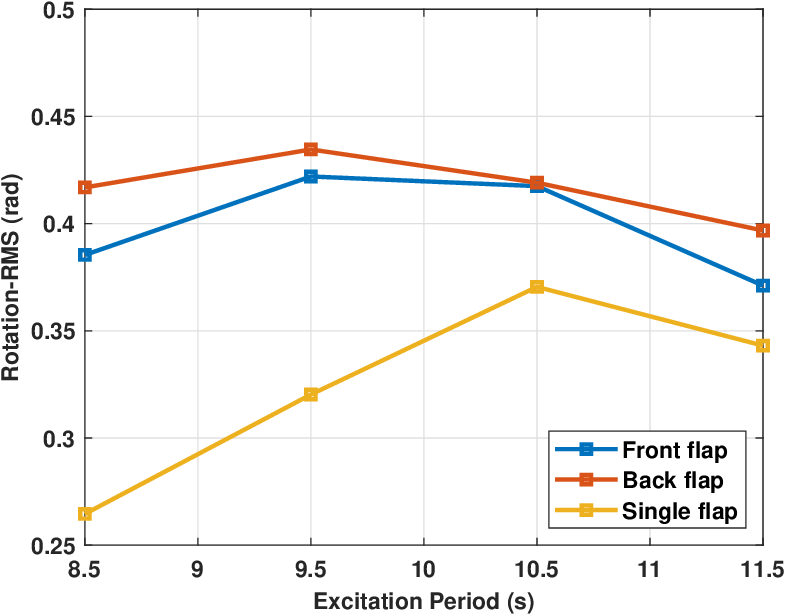}}%
\caption{The RMS values of the rotation of the flaps separated by (a) 10 m (b) 45m (c) 70 m under regular wave forcing with of 9.5s and wave height of 1.75m (left) and 3.25m (right)}
\label{fig: waves}
\end{figure}

The location of the WEC in the array can affect how much it is being impacted by the surrounding wave field. Wave-forced simulations allow for investigating how the front flap is impacted as the first flap to see the incident wave and how it is changing the field and the incident waves for the back flap when positioned with different separation distances. A series of simulations were conducted under different regular wave conditions with different wave periods ranging from 7.5s to 11.5s and two selected wave heights, 1.75 m and 3.25 m. Figure \ref{fig: waves} presents the impact of the separation distance on the flaps when excited with regular waves by comparing the rotation RMS of both flaps to the single one excited by the same wave conditions. The results show that for a separation distance of 10 m ($0.06\lambda < d < 0.11\lambda$), the front flap is more impacted by the wave field than the back one, and both of them are contributing to additional excitation from one on the other. For a higher amplitude, the flaps are experiencing destructive effects for periods longer than 10.5 s. For a separation distance of 45 m ($0.26\lambda < d < 0.51\lambda$), the back flap is starting to experience more impact from the wave field, especially at long periods. Finally, for a separation distance of 70 m ($0.41\lambda < d < 0.80\lambda$), the impact on both flaps is similar, and the change in their response is insignificant. It should be noticed that the difference in the rotation RMS of the flaps from the RMS of the single flap is decreasing as the separation distance increases. Overall, the flaps will have a constructive effect on each other under different regular wave conditions.

Comparing Figure \ref{fig: arbtr} with Figure \ref{fig: waves}, it can be noticed that the destructive effect from diffraction is compensated by the constructive effect from the incident waves. This can be expected since the wave amplitudes of the incident waves are much higher than those of diffracted and radiated waves. All the results from torque-excited and wave-excited simulations show that when the separation distance is within the range of $0.06\lambda < d < 0.11\lambda$, the responses of the flaps will be significantly impacted. It has a constructive impact in some cases and a destructive impact in others, but the total impact will still be constructive. Moreover, when the separation distance is in the range of $0.25\lambda < d < 0.80\lambda$, the flaps have a constructive effect on each other, and the differences in their rotation are decreasing and becoming insignificant with larger distances.

\subsection{Impact of the wave direction}

\begin{figure}[b!]
        \centering
        \subcaptionbox{}{\includegraphics[width=0.5\textwidth]{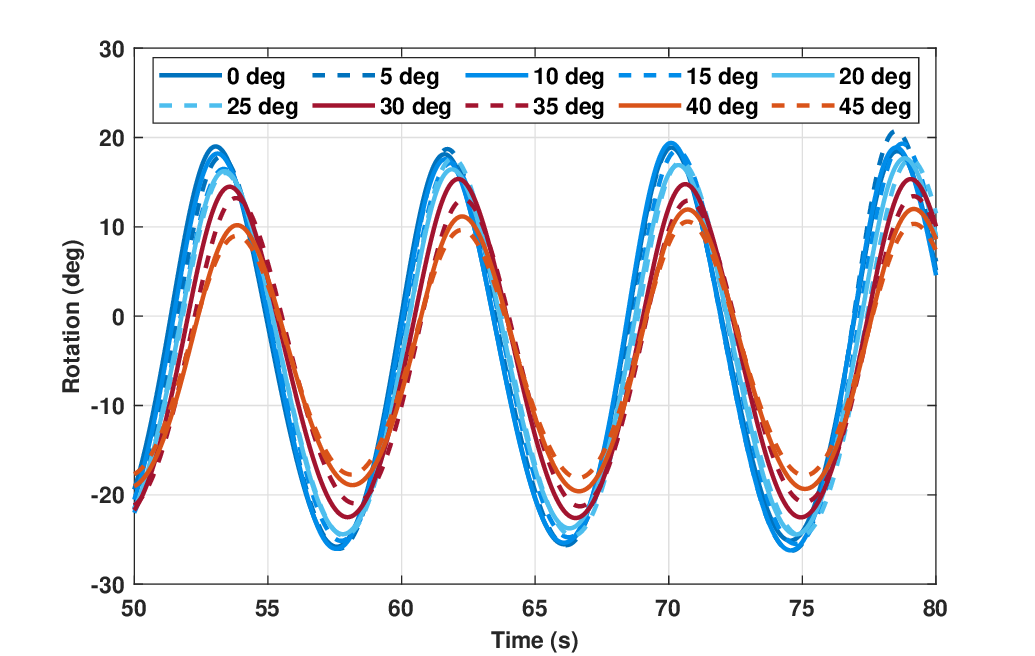}}%
        \hfill
        \subcaptionbox{}{\includegraphics[width=0.5\textwidth]{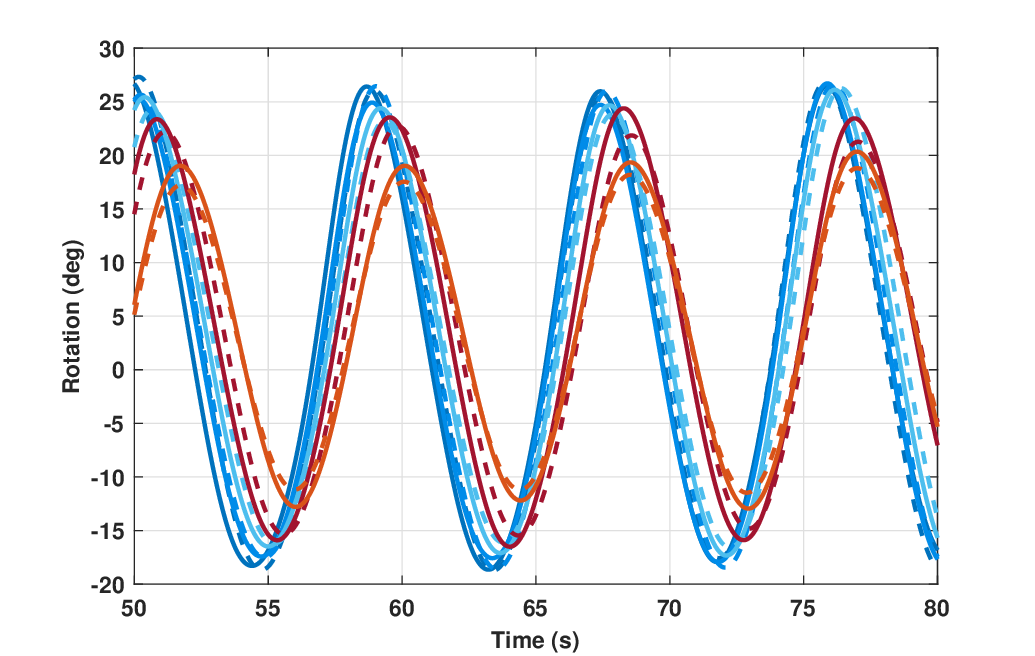}}%
         \caption{Comparison of simulated time series under regular wave forcing using a 8.5s wave period and 1.75m wave height obtained from different wave directions for (a) front and (b) back flaps with a separation distance of 45m between the flaps}
         \label{fig:waves angle}
\end{figure}

One of the main factors impacting the performance of WECs in arrays is the heading angle of the incident waves. For the flaps, the rotation amplitude is anticipated to be maximized when incident waves are perpendicular to their width and the heading angle is zero. However, changes in wave direction can alter flap rotation, potentially reducing its amplitude. To investigate the impact of the wave direction on the rotation amplitude of the flaps, a series of simulations was carried out with different heading angles for regular waves with an 8.5 s wave period and 1.75m wave height. The heading angles ranged from zero degree (perpendicular to flap width) to 45 degrees, incrementing a 5-degree step. This wave condition represents the most occurring wave at the targeted site, PacWave South, which means that it represents the greatest energy production loss.
Figure \ref{fig:waves angle} presents the rotation amplitude of the front and back flaps when separated by 45 m. The results show that when the heading angle is higher than 30 degrees, the rotation of the flaps is significantly impacted, resulting in an approximate 30\% loss in energy production. Notably, the back flap exhibits a slightly higher sensitivity to wave direction than the front flap, particularly in the phase angle. These findings come in agreement with the results in \cite{babarit2010impact}.
The study demonstrates the critical relationship between the heading angles of waves and the performance of WECs, highlighting a significant decrease in energy production as heading angles increase. These findings contribute valuable insights for optimizing WEC array configurations and improving overall system efficiency.

\subsection{Annual energy production}

All the changes impacting the response of the flaps will eventually impact the annual energy production of the system. To evaluate the impact of the separation distance on the annual energy production, different distances were tested under different regular wave conditions. Figure \ref{fig: PM} summarizes the estimated absorbed annual energy from each simulation with a real sea-state, while Table \ref{tab: Tot energy} shows the total annual absorbed energy for each of the considered seven separation distances. The results show a minimal effect on the produced energy when changing the distance between the flaps; this is in agreement with the results from the simulation. It confirms that for shorter distances, the constructive and destructive effects compensate for each other, keeping the dual-flap OSWEC in a constructive position when compared to a single flap. As the distance increases, the difference in the rotation of the flaps decreases, reaching about the same mean value. All these contribute to no change in the annual energy production with different separation distances. One can think that shorter distances are preferred because they can reduce manufacturing costs. However, other factors must be considered, such as the effect of the distance on the surge motion of the flaps and mooring loads. Although Table \ref{tab: Tot energy} implies that a separation distance of 33m is preferred as it generates the same amount of energy, there is a need for further investigation of the impact of using this distance on the surge motion and mooring loads. Our preliminary analysis shows a significant increase in both quantities and that the 45 m separation is a more optimized configuration when all factors are considered.

\begin{table}
\centering
    \caption{Annual mechanical energy from partial power matrix for different distances}
    \label{tab: Tot energy}
    \scalebox{1}{
\begin{tabular}{c c}
\hline
\textbf{Distance (m)} & \textbf{Mechanical Energy (GWh)} \\ \hline

Single (doubled)      & 1.15(2.3)             \\
10                    & 3.99                  \\ 
15                    & 3.76                  \\ 
33                    & 4.03                  \\ 
45                    & 3.98                  \\ 
55                    & 4.03                  \\ 
70                    & 4.07                  \\ 
86                    & 3.88                  \\ \hline
\end{tabular}}
\end{table}
  
\begin{figure}
\centering
\subcaptionbox{}{\includegraphics[clip, trim={3.5cm 0.3cm 0.0cm 1cm}, width=0.485\textwidth]{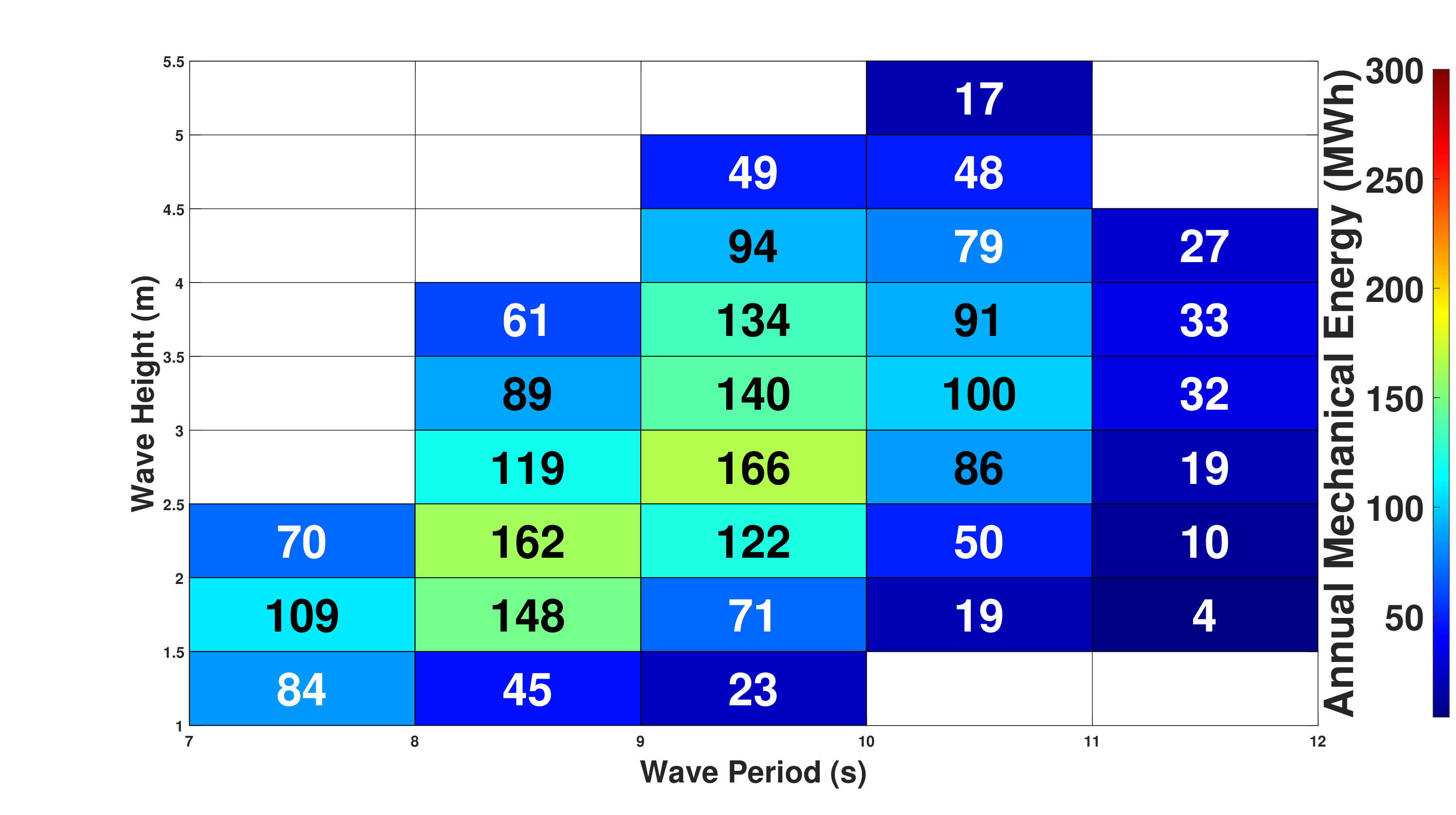}}%
\hfill
\subcaptionbox{}{\includegraphics[clip, trim={3.5cm 0.3cm 0.0cm 1cm}, width=0.485\textwidth]{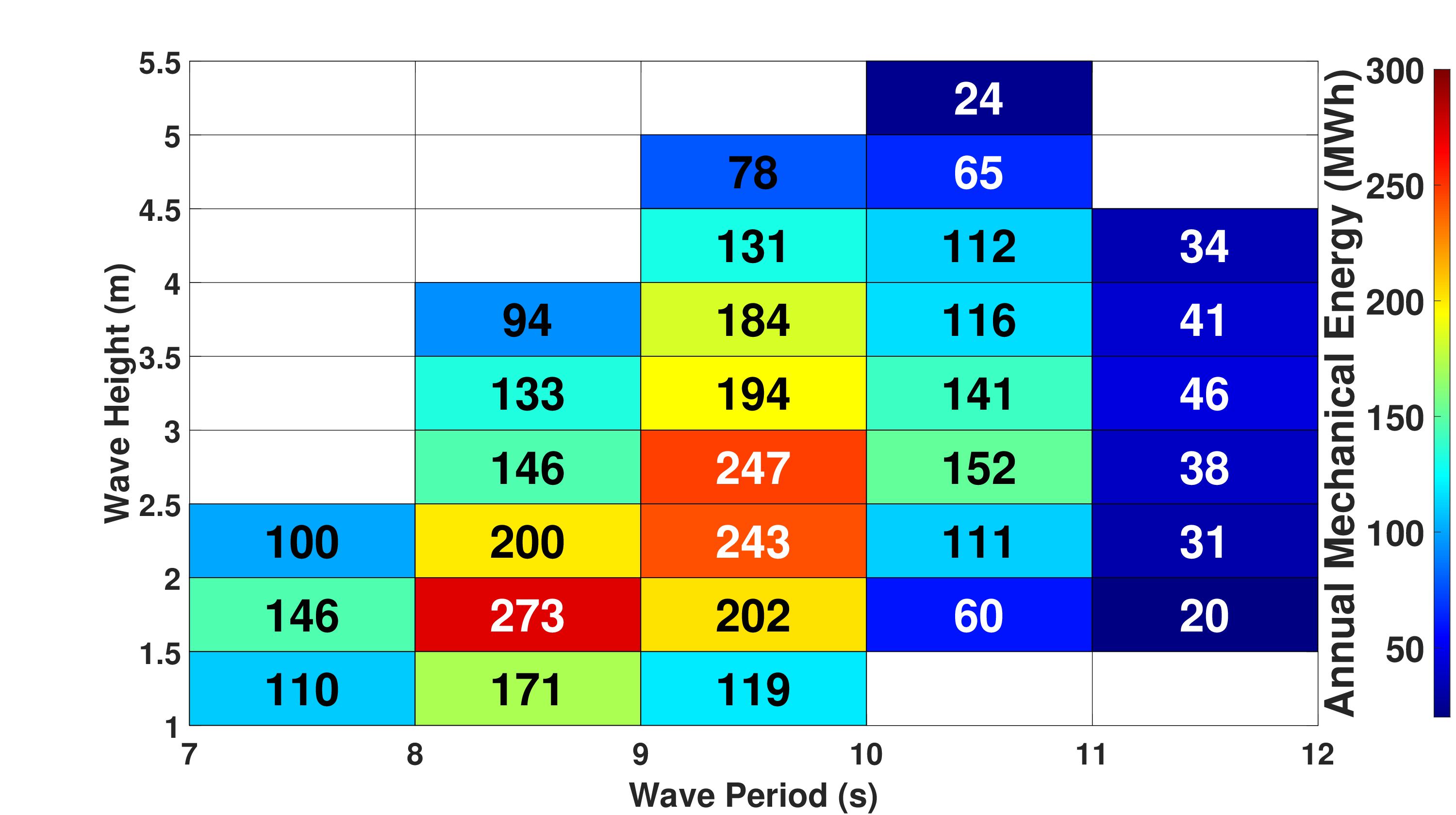}}%
\vfill
\subcaptionbox{}{\includegraphics[clip, trim={3.5cm 0.3cm 0.0cm 1cm}, width=0.485\textwidth]{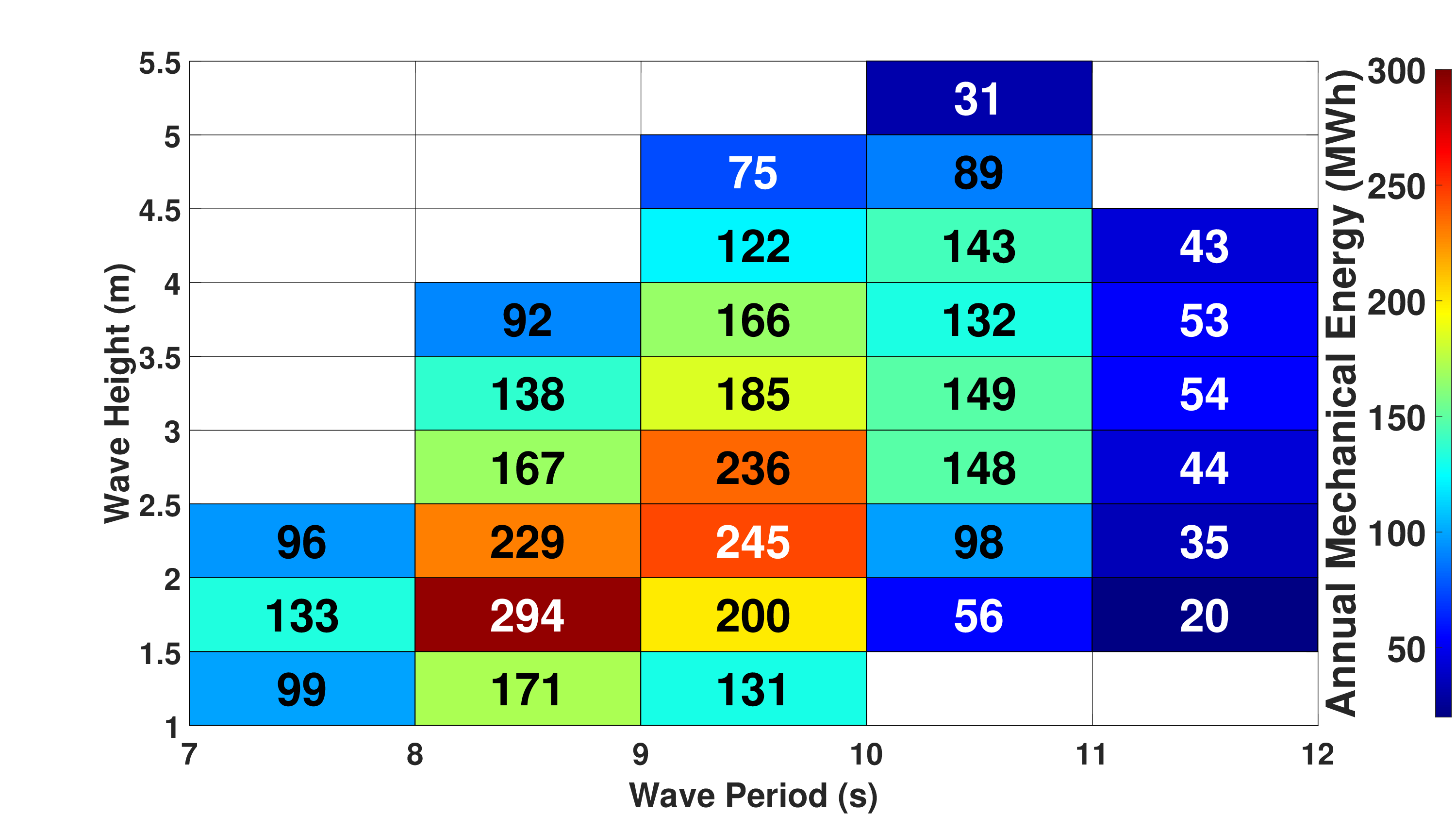}}%
\hfill
\subcaptionbox{}{\includegraphics[clip, trim={3.5cm 0.3cm 0.0cm 1cm}, width=0.485\textwidth]{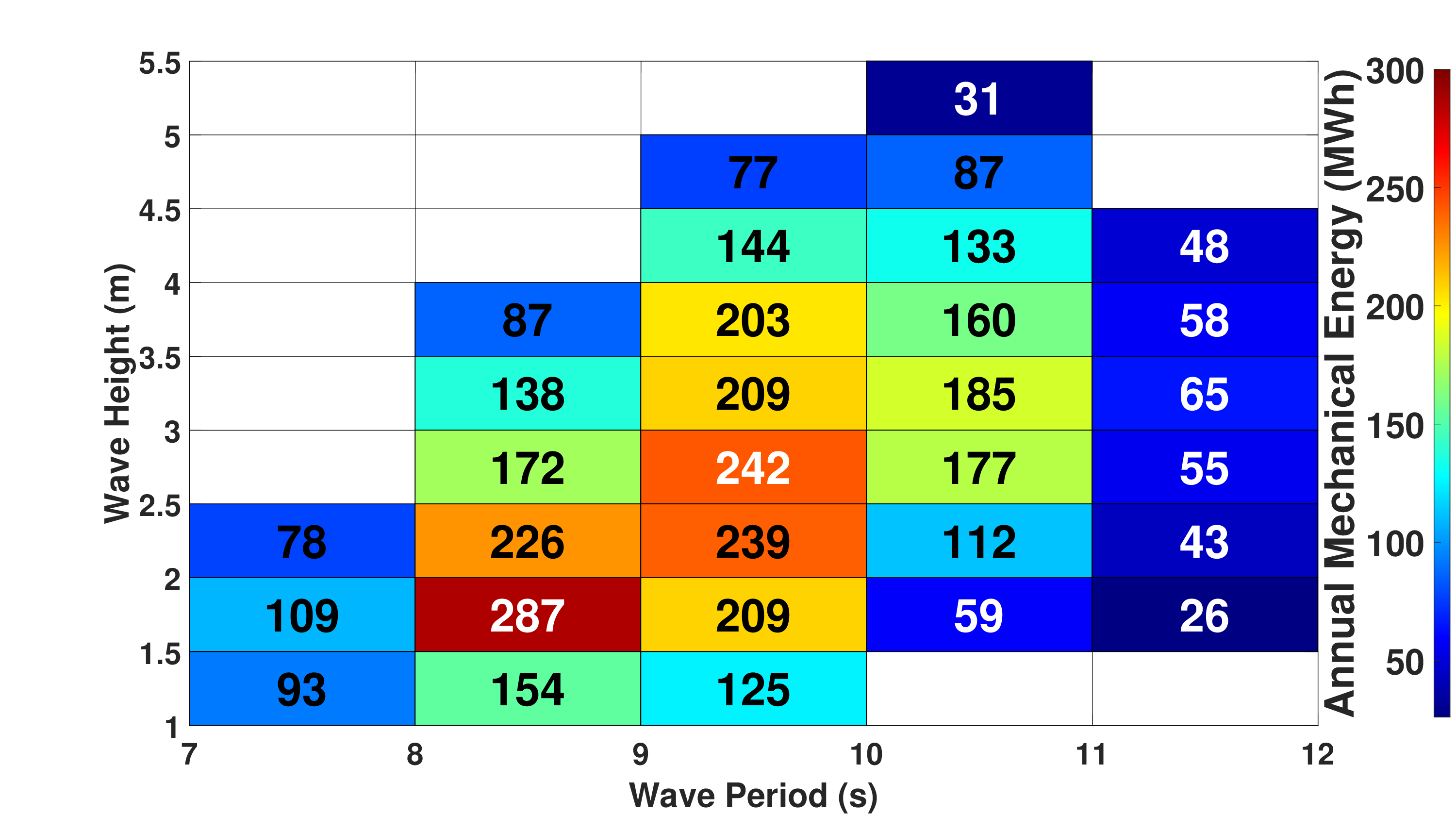}}%
\vfill
\end{figure}
\begin{figure}
\ContinuedFloat
\subcaptionbox{}{\includegraphics[clip, trim={3.5cm 0.3cm 0.0cm 1cm}, width=0.485\textwidth]{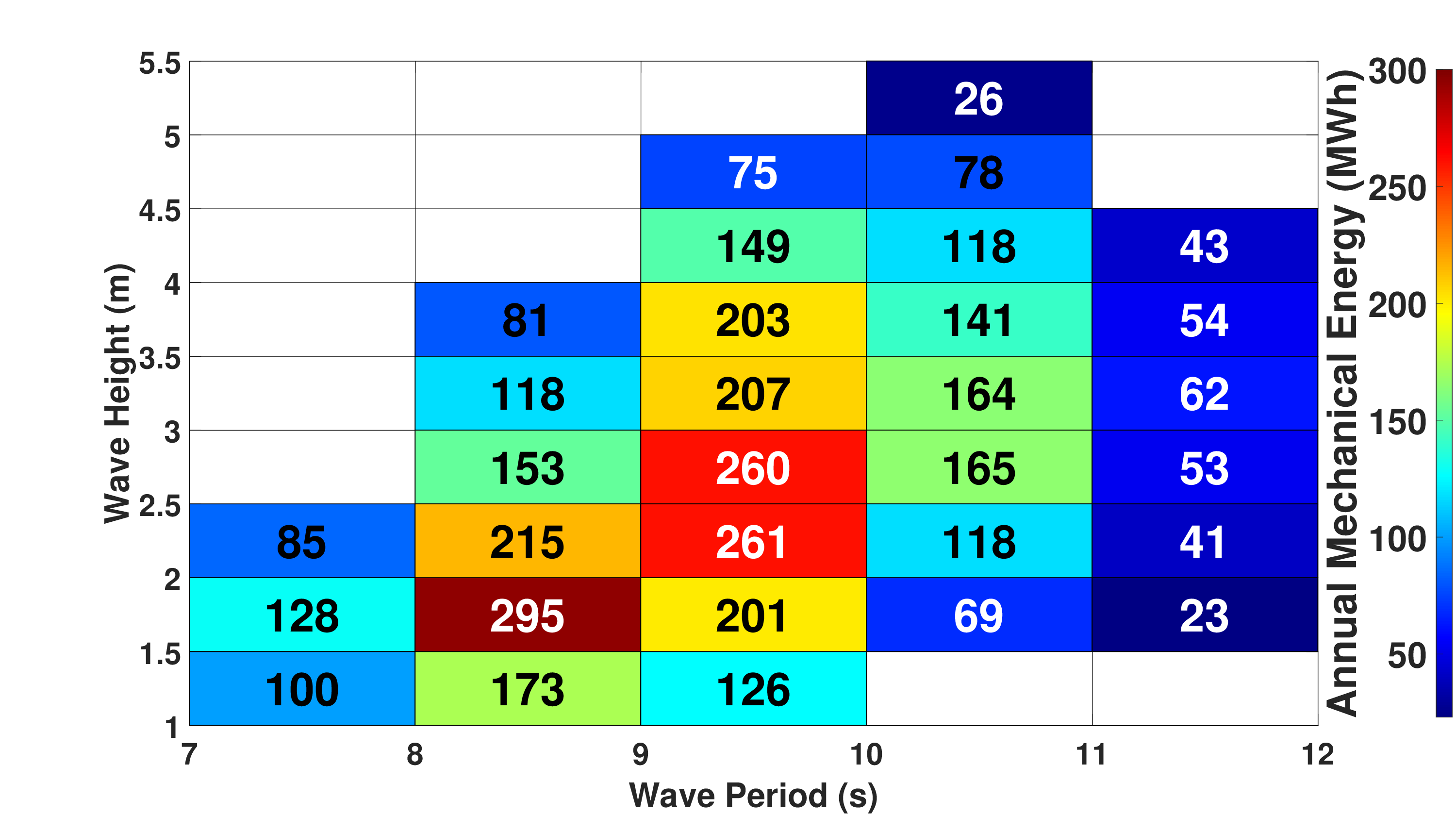}}%
\hfill
\subcaptionbox{}{\includegraphics[clip, trim={3.5cm 0.3cm 0.0cm 1cm}, width=0.485\textwidth]{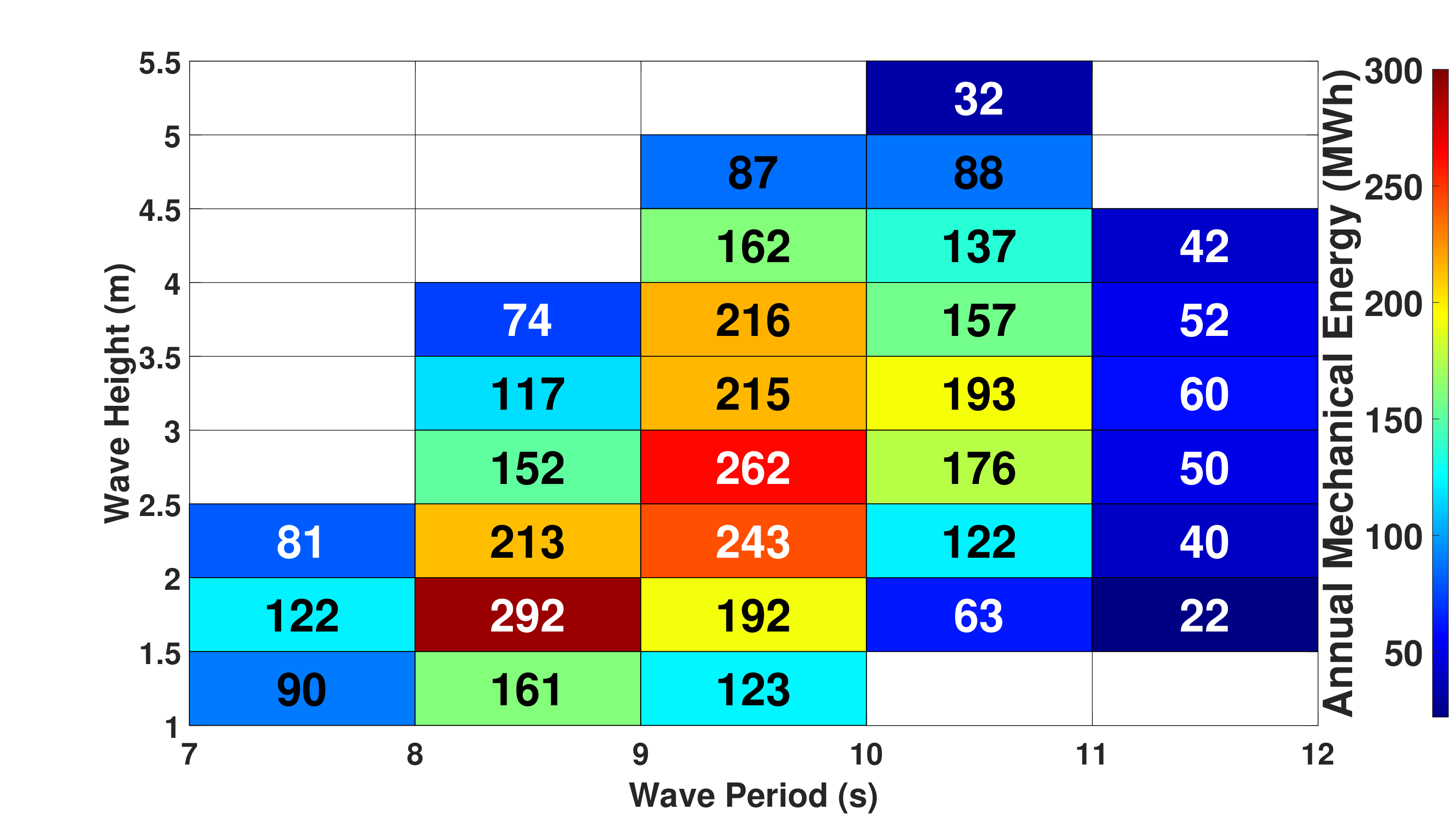}}%
\vfill

\subcaptionbox{}{\includegraphics[clip, trim={3.5cm 0.3cm 0.0cm 1cm}, width=0.485\textwidth]{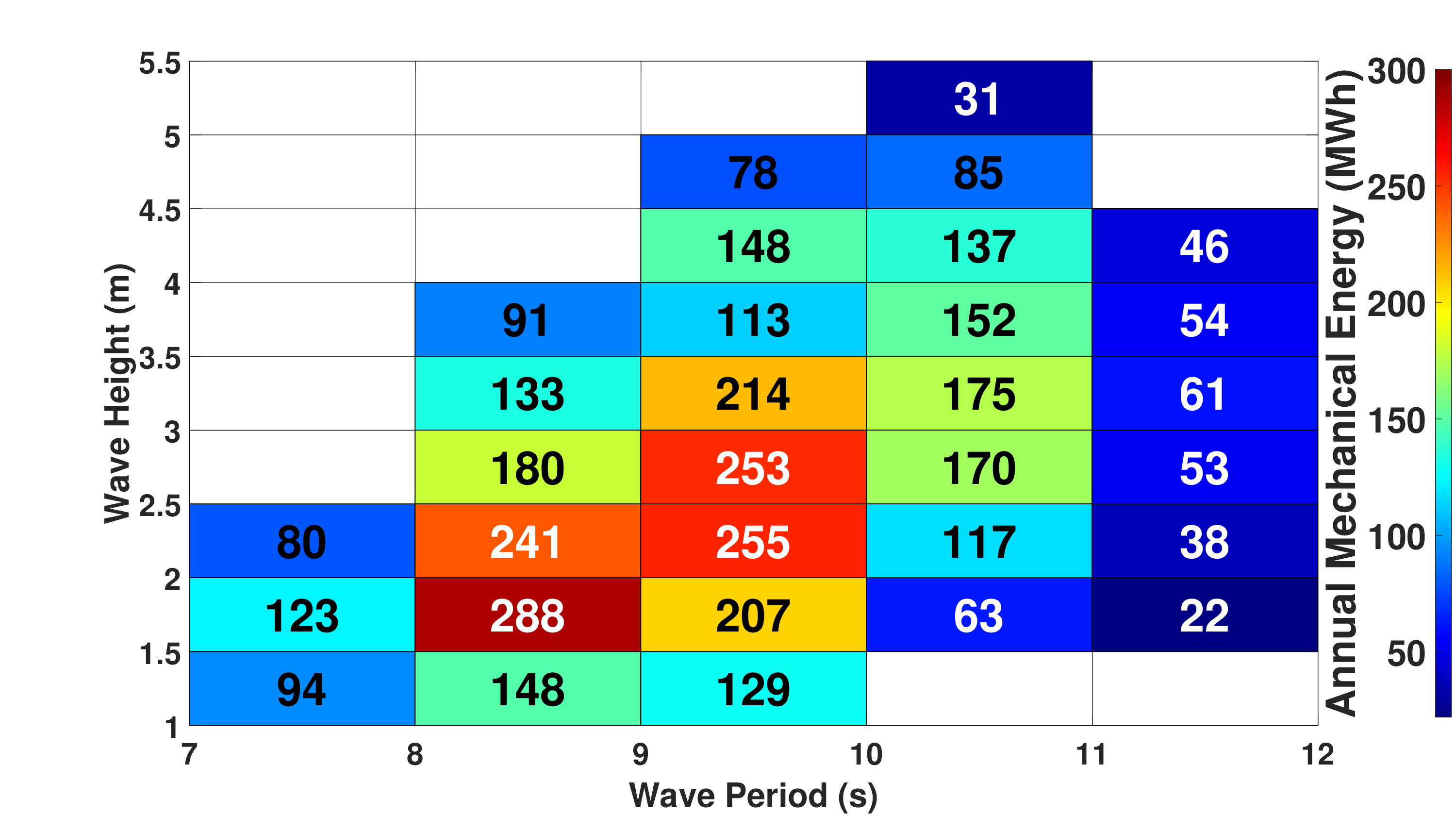}}%
\hfill
\subcaptionbox{}{\includegraphics[clip, trim={3.5cm 0.3cm 0.0cm 1cm}, width=0.485\textwidth]{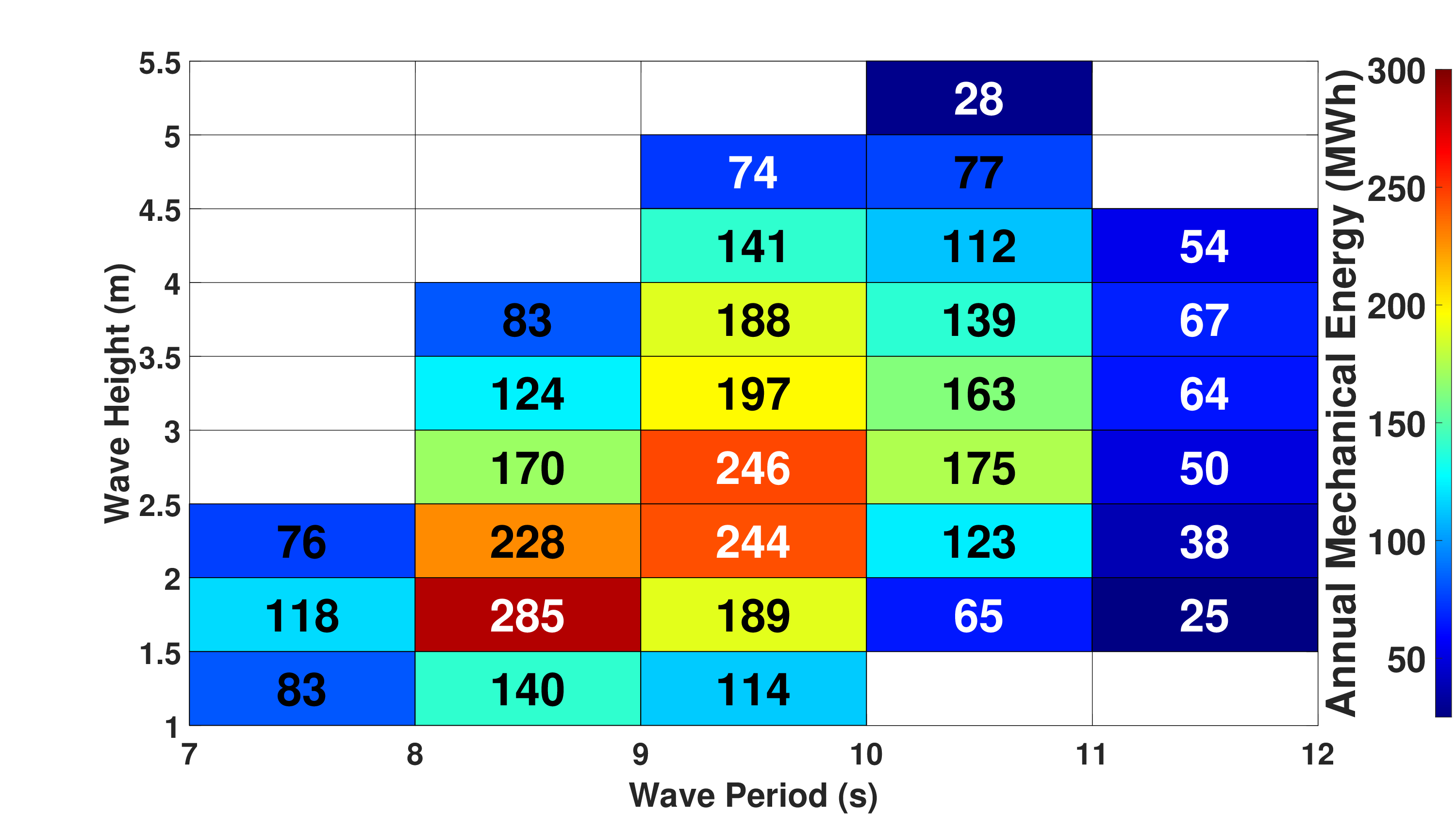}}%

\caption{Estimated annual mechanical energy of the (a) single flap multiplied by two and the dual-flap OSWEC with a separation distance of (b)10m (c)15m (d)33m (e)45m (f)55m (g)70m (h)86m for different wave conditions defined in JPD}
\label{fig: PM}
\end{figure}

\FloatBarrier
\section{Conclusions}

A dual-flap OSWEC was considered to investigate the impact of the separation distance between the flaps on the performance of each flap and the annual energy production. The system was treated as an array of two OSWECs; however, only short separation distances ranging from 10m to 86m were tested. The response of the flaps was numerically simulated under regular wave-forced and torque-forced conditions. To validate the numerical simulations, a 1:10 scaled model was manufactured and tested in the wave tank of the Davidson Laboratory at Stevens Institute of Technology. Torque-forced simulations were conducted to investigate the effect of each component in the equation of motion and its impact on the response of the flaps, while wave-forced simulations were conducted to investigate the effect of the location of the flap in the wave field and, consequently, its effect on the other flap. The results show that for a separation distance in the range of $0.06\lambda < d < 0.11\lambda$, the destructive impact from diffraction and the constructive impact from radiation are compensating each other. The interaction will be impacted by the phase difference of the response, causing constructive and destructive effects. However, the total impact will be constructive. The front flap is more impacted by the wave field from incident waves, but both flaps have constructive interference. For a separation distance in the range of $0.25\lambda < d < 0.80\lambda$, the destructive impact from diffraction is dominant over the constructive impact from radiation. The phase difference has a minimal effect on the response. Furthermore, incident waves have a constructive impact, compensating for that loss from diffraction and leaving both flaps with a constructive impact on each other. The change in their responses is decreasing, reaching the same mean value with an increasing separation distance. The impact of the wave direction was investigated. The energy production will be at its maximum with a zero heading angle and decrease to zero energy production when the waves become parallel to the device. The reduction in the response of the flaps is significant when the heading angle is greater than 30 degrees, resulting in an approximate 30\% loss in energy production. A partial power matrix was estimated for every separation distance based on the deployment of the dual-flap OSWEC at the targeted site, namely, the PacWave South site. The impact of the separation distance on the annual energy production is found to be insignificant, which gives more space for other factors to be satisfied in designing the arrays, such as mooring configuration. Further investigations are still needed to have a better insight on arrays, including but not limited to the wave direction and the layout of the array.

\section*{Acknowledgement}
This work was supported by the U.S. Department of Energy (DoE), Office of Energy Efficiency and Renewable Energy (EERE) under Award Number DE-EE0008953.\\
\indent Alaa Ahmed acknowledges the Link Foundation for the financial support through the Ph.D fellowship program in Ocean Engineering.

\bibliographystyle{cas-model2-names}

\bibliography{cas-refs}

\end{document}